\documentclass[twoside]{article}
\usepackage[english]{babel}
\usepackage{cite}
\usepackage{amssymb,amsmath,amsthm}
\usepackage{longtable}
\usepackage{nicefrac}
\usepackage[justification=justified,singlelinecheck=false]{caption}
\usepackage{color}
\usepackage[section]{placeins}
\newtheorem{defn}{Definition}[section]

\newtheorem{pro}{Proposition}[section]
\newtheorem{lem}{Lemma}[section]

\begin{document}
\begin{center}
{\Large \textbf{Isomorphism classes and invariants of low-dimensional filiform Leibniz
algebras}}\\[5mm]

{ A. O. Abdulkareem$^a$, I. S. Rakhimov$^b$ and S. K. Said Husain$^c$}

$ \ ^{a,b}$ Institute for Mathematical Research (INSPEM) and \\
$^{b,c}$Department of Mathematics, Faculty of Science, \\ Universiti Putra Malaysia \\ UPM, 43400,
Serdang, Selangor Darul Ehsan, Malaysia,\\[0pt]

$  ^{a}$abkaaz4success@yahoo.com \ $^{b}$risamiddin@gmail.com, \ $^c$skartini@science.upm.edu.my
\end{center}

\begin{abstract} In this study, we extend the result on classification of a subclass of filiform Leibniz algebras in low dimensions to dimensions seven and eight based on the technique used by Rakhimov and Bekbaev for classification of subclasses which arise from naturally graded non-Lie filiform Leibniz algebras. The subclass considered here arises from naturally graded filiform Lie algebras. This subclass contains the class of filiform Lie algebras and consequently, by classifying this subclass, we again re-examine the classification of filiform Lie algebras. Our resulting list of filifom Lie algebras is compared with that given by Ancoch\'ea-Bermudez and Goze in 1988 and by G\'{o}mez, Jimenez-Merchan and Khakimdjanov in 1998. 
\end{abstract}

\thispagestyle{empty}
\section{Introduction}
The concept of Leibniz algebra was introduced by Loday \cite{LP} in the study of Leibniz (co)homology as a noncommutative analogue of Lie algebras (co)homology. A Leibniz algebra over a field $K$ is a vector space over $K$ equipped with a $K$-bilinear map $[\cdot,\cdot]:L\times L\longrightarrow L$ satisfying the Leibniz identity  $[x,[y,z]]~=~[[x,y],z]~-~[[x,z],y],$  for
 all $ x,y,z\in L.$ Clearly, a Lie algebra is a Leibniz algebra, and conversely, a Leibniz algebra $L$ over the field $K$ ($char K \neq 2$) with property $[x,y]=-[y,x],$
for all $x,y\in L$ is a Lie algebra. Any Leibniz algebra $L$ gives rise to a Lie algebra $L_{Lie}$, which is obtained as the quotient of $L$ by the relation $[x,x]=0.$
 Let $I$ be the ideal of $L$ generated by all squares. Then $I$ is the minimal ideal with respect
to the property that $\mathfrak{g}:=L/I$ is a Lie algebra. The quotient map $\pi:L\longrightarrow \mathfrak{g}$ is a homomorphism of Leibniz algebras. 

The term filiform Lie algebra was introduced by M. Vergne \cite{MV} while studying the variety of nilpotent Lie algebras laws through which she intended to give a classification of nilpotent complex Lie algebras in dimension less than six (see also\cite{GJKh}).
The class of filiform Leibniz algebras in dimension $n$ have been split into three disjoint subclasses \cite{RB,GO}. In \cite{RB}, the subclasses have been denoted by $FLb_n$, $SLb_n$ and $TLb_n$ for the first, second and third classes respectively.

The first and second classes which arise from naturally graded non-Lie filiform Leibniz algebra have been classified up to dimension $ n \leq 9,$ \cite{Fatanah, SIK} and the references therein. 

The present paper is concerned with a subclass of filiform Leibniz algebras. This subclass arises from naturally graded filiform Lie algebras. In dimension $n$, it has been denoted by $TLb_n$ and treated before in \cite{OR, RH1, RH2}. In case of $TLb_n$, the ideal $I$ generated by squares coincides with the right center of the Leibniz algebras since we are considering right Leibniz algebras. On the contrary, the left center is only a subalgebra in the context of right Leibniz algebras. Therefore, derived algebras are one-dimensional Leibniz central extensions of the Lie algebra $L_{Lie}.$  The aim of this paper is to give complete lists of isomorphism classes of one-dimensional Leibniz central extensions of filiform Lie algebras of dimensions $6$ and $7$. As for the lower dimensions they either have been treated before \cite{RH1} or have been included in the list of nilpotent Leibniz algebras \cite{AOR1}, \cite{AOR2}. Another motivation to study $TLb_n$ comes out from the fact that it contains the class of all $n$-dimensional filiform Lie algebras. Classifying $TLb_n$, we once again examine thoroughly and recover the lists of low-dimensional filiform Lie algebras given in \cite{AG,GJKh}.

Organization of this paper is as follows. Section 2 is devoted to providing necessary results needed throughout the study. In the third and fourth sections, complete classification of $TLb_7$ is given with the proof of some propositions that are stated therein. Section $5$ and $6$ are devoted to list the isomorphism classes of $TLb_8$. We describe the orbits of the base change for single orbits and for union of parametric family of orbits cases. In parametric cases, invariant functions are given to distinguish the orbits. Necessary conclusions drawn from this study form the contents of Section $7$.
\section{Preliminaries}
\label{sec:2}
In this section we give some facts on Leibniz algebras that will be
used throughout the paper.

\begin{defn}
An algebra $L$ over a field $K$ is called a Leibniz algebra, if its
bilinear operation $\lambda(\cdot,\cdot)$ satisfies the following Leibniz
identity:
\begin{equation}
\lambda(x,\lambda(y,z))=\lambda(\lambda(x,y),z)-\lambda(\lambda(x,z),y), \ for \  x,y,z\in L.
\end{equation}
\end{defn}
A Leibniz algebra structure on $n$-dimensional vector space $V$ over a field $K$ can be regarded
as a pair $L=(V,\lambda),$ where $\lambda$ is a Leibniz algebra law on $V$. We denote by $LB_n$ the set of all Leibniz algebra laws on $V$. It is a subspace of the linear space of $Hom(V\otimes V, V)$.
\begin{defn}  An action of a group $G$ on a set $X$ is a function $\ast:G\times X\rightarrow X$ that satisfies the following conditions;
\begin{enumerate}
\item $e\ast x=x$, $\ for \ all \ x \in X$, where $e$ is the identity element of $G$.
\item $g\ast(h\ast x))=(gh)\ast x$, for all $g,h \in G$ and $x\in X$.
\end{enumerate}
\end{defn}
The above group action is a morphism of algebraic variety when $G$ and $X$ are substituted for algebraic group and algebraic variety respectively. In what follows, we define $(g,\lambda)\mapsto g\ast \lambda$

\begin{defn} Two laws $\lambda_1$ and $\lambda_2$ from $LB_n$ are said to be isomorphic, if there exists $g\in GL(V)$ 
such that
$$\lambda_2(x,y) = (g\ast \lambda_1)(x,y) = g^{-1}(\lambda_1(g(x),g(y )))$$ for all $x, y \in V$. 
\end{defn}
Thus we get an action of $GL_n(K)$ on $LB_n$.
Let $O(\lambda)$ be the set of the laws isomorphic to $\lambda,$ it is called the orbit of $\lambda.$
Let us fix a basis $\{e_1, e_2, e_3,..., e_n\}$ of $V.$ The structure constants, $\gamma_{ij}^k \in K,$ of $\lambda \in LB_n$ are given by
$$\lambda(e_i, e_j) =\sum \limits_{k=1}^n \gamma_{ij}^k e_k, \ \ i,j=1,2,...,n.$$
Once a basis is fixed, we can identify the law $\lambda$ with its structure constants. These
constants are the structure constants of a Leibniz algebra if and only if:
$$\sum\limits_{\emph{l}=1}^{\emph{n}}{(\gamma_{\emph{jk}}^{\emph{l}}
\gamma_{\emph{il}}^{\emph{m}}-\gamma_{\emph{ij}}^{\emph{l}}\gamma_{\emph{lk}}^{\emph{m}}+
\gamma_{\emph{ik}}^{\emph{l}}\gamma_{\emph{lj}}^{\emph{m}})}=0, \ \ i,j,k,m=1,2,...,n.$$
Then $LB_n$ appears as a subvariety of the algebraic variety,\, $Alg_n(K)$,\, of all algebraic structures on  $V$.
Let $\lambda\in LB_n$ and $G_\lambda$ be the subgroup of $GL_n(K)$ defined by
$G_\lambda~=~\{g~\in~GL_n(K)~|~g~\ast~\lambda~=~\lambda\}.$
The orbit $O(\lambda)$ of $\lambda$ with respect to the action of $GL_n(K)$  is isomorphic to the homogeneous space
$GL_n(K)/G_\lambda.$
\begin{defn}
A function $f:LB_n\longrightarrow K$ is said to be invariant (or orbit) function with respect to an action of a subgroup $G$ of $GL_n(K)$ on $LB_n$ if $f(g\ast\lambda)=f(\lambda)$ for $g\in G.$
\end{defn}
In the foregoing, all algebras are assumed to be over the fields of complex numbers $\mathbb{C},$ and the bracket
notation is used to describe the law: $[x,y]=\lambda(x,y)$ for $x,y\in L.$

Let $L$ be a Leibniz algebra. Define
\begin{equation}
L^{1}=L,\ L^{k+1}=[L^{k},L],\ k\geq 1. 
\end{equation}
Clearly, $$ L^{1}\supseteq L^{2}\supseteq \cdots $$
\begin{defn}
A Leibniz algebra L is said to be a nilpotent, if there exists $s\in
\mathbb{N},$ such that
\begin{equation}
L^{1}\supset  L^{2}\supset ...\supset L^{s}=\{0\}.
\end{equation}
\end{defn}
\begin{defn}
A Leibniz algebra $L$ is said to be a filiform, if $dimL^{i}=n-i,$ where $%
n=dimL$ and $2\leq i\leq n.$
\end{defn}

According to \cite{OR,AOR2} the class of $(n+1)$-dimensional filiform
Leibniz algebras arising from naturally graded filiform Lie
algebras, denoted by $TLb_{n+1}$, admits a so-called adapted basis $\{e_0,e_1,...,e_n\}$
such that the table of multiplication with respect to the basis has the following form:
 \begin{eqnarray} \label{TLb} TLb_{n+1}=\left\{
\begin{array}{lll}
\lbrack e_{i},e_{0}]=e_{i+1}, \qquad \qquad \qquad \qquad \qquad \qquad \qquad \qquad \quad 1\leq i\leq {n-1},  \\[1mm]
\lbrack e_{0},e_{i}]=-e_{i+1}, \qquad \qquad \qquad \qquad \qquad \qquad \qquad \qquad  \   2\leq i\leq {n-1}, \\[1mm]
\lbrack e_{0},e_{0}]=b_{0,0}e_{n},  \\[1mm]
\lbrack e_{0},e_{1}]=-e_{2}+b_{0,1}e_{n}, \\[1mm]
\lbrack e_{1},e_{1}]=b_{1,1}e_{n},  \\[2mm]
\lbrack e_{i},e_{j}]={a_{ij}^{1}e_{i+j+1}+\dots
+a_{ij}^{n-(i+j+1)}e_{n-1}}+b_{ij}e_{n}, \quad  1\leq i<j\leq {n-1},   \\[2mm]
\lbrack e_{i},e_{j}]=-[e_{j},e_{i}],\quad \qquad \qquad \qquad \qquad \qquad \qquad \qquad \, 1\leq i<j\leq n-1,  \\[1mm]
\lbrack
e_{i},e_{n-i}]=-[e_{n-i},e_{i}]=(-1)^{i}b_{i,n-i} e_{n}, \\[1mm]
 where \; a_{ij}^k, b_{ij} \in \mathbb{C}, and $\ $b_{i,n-i}=0, \quad \qquad \qquad \qquad \quad \; 1\leq i\leq {n-1},\\  for \ even $ \ $n.
\end{array}%
\right.\end{eqnarray}
In the table (\ref{TLb}), there are interrelations among $a_{ij}^k, b_{ij}$ that are different in each fixed dimension as evident in the dimensions considered here.

The class $TLb_{n+1}$ can be regarded as follows. Let $\mu_{n}$ be the $n$-dimensional Lie algebra defined by the
brackets: $[e_i,e_0]=e_{i+1}, \ i=1,2,...,n-2$ and
$\{e_0,e_1,...,e_{n-1}\}$ be a basis of $\mu_{n}$. Clearly, $\mu_{n}$ is a filiform Lie
algebra. Moreover, any $n$-dimensional filiform Lie algebra is
isomorphic to a linear deformation of $\mu_{n}$ (see \cite{GJKh}).
More precisely, let $\Delta$ be the set of all pairs of integers
$(k,r)$ such that $1\leq k \leq \left[\frac{n-2}{2}\right],\ \
2k+1<r\leq n-1,$ (if $n$ is odd then $\Delta$ contains also the
pair $\left(\frac{n-2}{2},n-1\right)$). For any elements $(k,r)\in
\Delta,$ one can associate the 2-cocycle, $\Psi_{k,r}$ for the Chevalley
cohomology of $\mu_n$ with coefficients in the adjoint module defined by 
 \begin{equation}
\Psi_{k,r}(e_i,e_j)=\left\{\begin{array}{ll}-\Psi_{k,r}(e_j,e_i)=(-1)^{k-i} {j-k-1\choose
k-i}  \,e_{i+j+r-2k-1},\ 1\leq i \leq k <j \leq n-1 & \\ \ \ \ 0, \quad \mbox{otherwise} 
\end{array}\right.
\end{equation}
Any $n$-dimensional filiform Lie algebra law $\mu$ in the variety of
filiform Lie algebras laws is isomorphic to $\mu_n+\Psi,$ where $\Psi$ is a 2-cocycle defined by
$$\Psi=\sum_{\left(k,r\right)\in \Delta}a_{k,r}\,\Psi_{k,r},$$ and verifying
the relation $\Psi \circ\Psi=0$ with
\begin{equation}
(\Psi \circ\Psi)(x,y,z)= \Psi(\Psi(x,y),z)+\Psi(\Psi(y,z),x)+\Psi(\Psi(z,x),y).
\end{equation}
In this situation, we consider another cocycle which is Leibniz 2-cocycle, $\Psi$ satisfying the relation,
\begin{equation}
(\Psi \circ\Psi)(x,y,z)=
\Psi(x,\Psi(y,z))+\Psi(\Psi(x,z),y)-\Psi(\Psi(x,y),z)=0.
\end{equation}
One-dimensional Leibniz central extension of linear deformation of the Lie algebra $\mu_n$ has been denoted by {\em Ced($\mu_n$)} in \cite{RH3}, the resulting table of multiplication has been slightly modified and is given by 
\begin{equation}
 Ced(\mu_{n})=\left\{
 \begin{array}{lll}
 \lbrack e_{i},e_{0}]=e_{i+1},\quad \qquad \qquad \qquad \qquad \qquad \qquad \qquad \qquad  \qquad \qquad  \ 1\leq i\leq {n-1}, \\
 \lbrack e_{0},e_{i}]=-e_{i+1},\; \qquad \qquad \qquad  \qquad \qquad \qquad \qquad \qquad \qquad \qquad \ 2\leq i\leq {n-1}, \\
 \lbrack e_{0},e_{0}]=b_{00}e_{n}, \\ \lbrack e_{0},e_{1}]=-e_{2}+b_{01}e_{n},  \\
 \lbrack e_{1},e_{1}]=b_{11}e_{n},  \\
 \lbrack e_{1},e_{i}]=-[e_i,e_1]=a_{1,i+2}e_{i+2}+a_{1,i+3}e_{i+3}+\cdots+b_{1,i}e_n,\ \qquad \; 2\leq i\leq n-3,\\
 \lbrack e_{i},e_{j}]=-[e_j,e_i]\in span\{e_{i+j+1},e_{i+j+2},\cdots,e_{n}| \ 2\leq j\leq n-j-1\} \\
 \lbrack e_{i},e_{n-i}]=-[e_{n-i},e_{i}]=(-1)^{i}b_{i,n-i}e_{n}, \\
 \mbox{ where } a_{i,j}^{k},b_{i,j}\in \mathbb{C} \mbox{ and } b_{i,n-i}=b \mbox{ whenever }\ \qquad \qquad  \qquad \qquad  1\leq i \leq n-1,\\ \mbox{ and }  b=0 \ for \mbox{ even } n.
 \end{array}%
 \right. 
 \end{equation}
\begin{defn}
Let $\{e_0,e_1,...,e_n\}$ be an adapted basis of $L\in
TLb_{n+1}.$ Then a nonsingular linear transformation
$f:L\rightarrow L$ is said to be adapted if the basis
 $ \{f(e_0),f(e_1),...,f(e_n)\} $ is adapted.
\end{defn}

The subgroup of $GL_{n+1}$ consisting of all adapted
transformations, we denote by $G_{ad}.$ The following proposition proved in\,\cite{RH1}
specifies the elements of $G_{ad}.$
\begin{pro}
Let $\{e_0,e_1,...,e_n\}$ be an adapted basis of $L\in TLb_{n+1}.$ and $f$ be an adapted transformation. Then $f$ can be represented as follows:
\begin{eqnarray}
f(e_0)&=& e_{0}^{\prime }=\sum_{i=0}^{n}A_{i}e_{i},\nonumber \\
f(e_1)&=& e_{1}^{\prime }=\sum_{i=1}^{n}B_{i}e_{i},\nonumber \\
f(e_i)&=& e_{i}^{\prime }=[f(e_{i-1}),f(e_0)], \qquad
2\leq i \leq n \nonumber
\end{eqnarray}
 $A_0, A_i, B_j,\ (i,j=1,...,n)$ {\em are complex
numbers and} $A_0\,B_1(A_0+A_1b)\neq 0.$
\end{pro}
In $G_{ad}$ we specify the following transformations called
elementary:
\begin{eqnarray}
\tau( a,b,c)&=&\left\{\begin{array}{lll} \tau(e_0)=a\,e_0+b\,e_1, &
\\[1mm]
\tau(e_1)=c\,e_1,\quad a\,c\neq 0,\\[1mm]
\tau(e_{i+1})=[\tau(e_i), \tau(e_0)], \quad 1\leq i\leq n-1,\nonumber  \\[2mm]
\end{array} \right.\\[2mm]
\sigma(a,k)&=&\left\{\begin{array}{lll} \sigma(e_0)=e_0+a\,e_k,
\qquad \ \ \quad 2\leq k \leq n, &
\\[1mm]
\sigma(e_1)=e_1, & \\[2mm]
\sigma(e_{i+1})=[\sigma(e_i), \sigma(e_0)], \quad 1\leq i\leq n-1,\nonumber  \\[2mm]
\end{array} \right.\\[2mm]
\phi(c,k)&=&\left\{\begin{array}{lll}\phi(e_0)=e_0, &
\\[2mm]
\phi(e_1)=e_1+c\,e_k, \quad \qquad \ \ \ 2\leq k \leq n,
\\[2mm]
\phi(e_{i+1})=[\phi(e_i), \phi(e_0)], \quad 1\leq i\leq n-1,\nonumber
\end{array}\right.
\end{eqnarray}
where $a,b,c \in \mathbb{C}$.

\begin{pro}\label{p}  
Let $L$ be an algebra from ${{TLb}}_{n+1}$, then any adapted transformation $f$ of $L$ can be represented
as the composition:
\begin{eqnarray} 
f=&\phi(B_n,n)\circ\phi(B_{n-1},n-1)\circ...\circ\phi(B_2,2)\circ\sigma(A_n,n)\nonumber\\& \circ \sigma(A_{n-1},n-1) \circ...\circ\sigma(A_2,2)\circ\tau(A_0,A_1,B_1)
\end{eqnarray}
\begin{proof} 
The proof can be obtained by direct verification.
\end{proof}
\end{pro}

 \begin{pro}\label{p1} The transformation
\begin{eqnarray}
g=&\phi(B_n,n)\circ\phi(B_{n-1},n-1)\circ\phi(B_{n-2},n-2)\circ... \circ\sigma (A_n,n)\nonumber\\&\circ\sigma(A_{n-1},n-1)\circ...\circ\sigma(A_{2},{2})
\end{eqnarray}
if $n$ even,  and
\begin{eqnarray}
g=&\phi(B_n,n)\circ\phi(B_{n-1},n-1)\circ\sigma(A_n,n)\circ\sigma(A_{n-1},n-1)\nonumber\\&\circ \sigma(A_{n-2},{n-2})\circ...\circ\sigma(A_{2},{2})
\end{eqnarray}
 for \ odd \ $n$, do not change the structure constants of algebras from $TLb_{n+1}.$
\end{pro}

For the proof of Proposition \ref{p1} we refer the reader to
\cite{RH2}. The proof of the following lemma can be found easily by a simple computation.
\begin{lem}\label{L1} Let $\{e_0,e_1,...,e_n\}\longrightarrow \{e_0',e_1',...,e_n'\}$ be an adapted base change, $b_{00},b_{01},b_{11},...$ and $b'_{00},b'_{01},b'_{11},...$ be the respective structure constants. Then for
$b'_{00},b'_{01}$\ and \ $b'_{11}$ one has
\begin{eqnarray}
b_{00}^\prime&=&\frac{A_0^2b_{00}+A_0A_1b_{01}+A_1^2b_{11}}{A_0^{n-2}B_1(A_0+A_1\,b)},\ \nonumber
 \\[2mm]
 b_{01}^\prime&=&\frac{A_0b_{01}+2A_1b_{11}}{A_0^{n-2}(A_0+A_1\,b)},\
\nonumber \\ b_{11}^\prime&=&\frac{B_1b_{11}}{A_0^{n-2}(A_0+A_1\,b)}.\nonumber
\end{eqnarray}
\end{lem}
\section{Isomorphism criterion for ${TLb}_{7}$}
Due to the observation in \cite{RH3} as pointed out in Section \ref{sec:2}, an algebra $L$ from ${TLb}_{7}$ is represented as a one dimensional Leibniz central extension\,($C(L)=<e_6>$)
 of the filiform Lie algebra with law $\mu_6+\Psi,$ where $\Psi$ is the appropriate cocycle on the adapted basis $\{e_0,e_1,...,e_6\},$ the class ${TLb}_{7}$ is represented as follows:
 \begin{center}
 $TLb_{7}=\left\{
 \begin{array}{lll}
 [e_{i},e_{0}]=e_{i+1}, \qquad \ \ \  1\leq i\leq
 {5},   \\[0mm]
 [e_{0},e_{i}]=-e_{i+1}, \qquad 2\leq i\leq {5},   \\[0mm]
 [ e_{0},e_{0}]=b_{00}e_{6},  \\[0mm]
 [ e_{0},e_{1}]=-e_{2}+b_{01}e_{6}, \\[0mm]
 [ e_{1},e_{1}]=b _{11}e_{6}, \\[0mm]
 [ e_{1},e_{2}]=-[e_2,e_1]=a_{14}\,e_4+a_{15}\,e_5+b_{12}e_{6} , \\[0mm]
 [ e_{1},e_{3}]=-[e_{3},e_{1}]=a_{14}\,e_5+b_{13}e_{6}, \\[0mm]
  [e_{1},e_{4}]=-[e_{4},e_{1}]=-a_{25}\,e_5+b_{14}e_{6}, \\[0mm]
  [e_{2},e_{3}]=-[e_{3},e_{2}]=a_{25}\,e_5+b_{23}e_{6},\\[0mm]
  [e_{1},e_{5}]=-[e_{5},e_{1}]=b_{15}e_{6}, \\[0mm]
  [e_{2},e_{4}]=-[e_{4},e_{2}]=b_{24}e_{6}.\\[0mm]
 \end{array}%
 \right. $ 
 \end{center}
Observe that the structure constants $a_{ij}$ and $b_{st}$ in the table above are not free. The interrelations between them are defined as follows.
\begin{lem}
 The structure constants of algebras from ${TLb}_{7}$ satisfy the following constraints: 
 \begin{enumerate}
 \item[{\em 1.}] $b_{13}=a_{15}$, 
 \item[{\em 2.}] $b_{14}=a_{14}-b_{23}$ \mbox{and}
 \item[{\em 3.}] $b_{15}=b_{24}=a_{25}=0$.
\end{enumerate}


\begin{proof}
The first and  second  relations occur if we apply the Leibniz identity to
the triples $\{e_0,e_1,e_2\},\
\{e_0,e_1,e_3\}$ of the basis vectors  respectively, and
the last relation is a result of application of the Leibniz identity to $\{e_0,e_2,e_3\},\{e_0,e_1,e_4\},$ and $\{e_1,e_2,e_3\}.$
\end{proof}
\end{lem}
Further unifying the table of multiplication above we rewrite it
via new notations
$c_{00},c_{01},c_{11},c_{12},c_{13},c_{14},c_{23}$ for the structure constants:
\begin{center}
$TLb_{7}=\left\{
\begin{array}{lll}
\lbrack e_{i},e_{0}]=e_{i+1},\quad \qquad 1\leq i\leq {5}, &  \\
\lbrack e_{0},e_{i}]=-e_{i+1},\; \qquad 2\leq i\leq {5}, &  \\
\lbrack e_{0},e_{0}]=c_{00}e_{6}, &  &  \\
\lbrack e_{0},e_{1}]=-e_{2}+c_{01}e_{6}, &  &  \\
\lbrack e_{1},e_{1}]=c _{11}e_{6}, &  &  \\
\lbrack e_{1},e_{2}]=-[e_2,e_1]=c_{12}\,e_4+c_{13}\,e_5+c_{14}e_{6}, \\
\lbrack e_{1},e_{3}]=-[e_{3},e_{1}]=c_{12}\,e_5+c_{13}e_{6}, \\[0mm]
 [e_{1},e_{4}]=-[e_{4},e_{1}]=(c_{12}-c_{23})\,e_{6}, \\[0mm]
 [e_{2},e_{3}]=-[e_{3},e_{2}]=c_{23}\,e_{6}.
\end{array}%
\right. $ 
 \end{center}
An algebra from ${TLb}_{7}$ with the structure constants
$c_{00},c_{01},c_{11},c_{12},c_{13},c_{14},c_{23}$ is
denoted by $L(C),$ where
$C=(c_{00},c_{01},c_{11},c_{12},c_{13},c_{14},c_{23}).$
The next theorem represents the action of the adapted base change on the structure constants
$c_{00},c_{01},c_{11},c_{12},c_{13},c_{14},c_{23}$ of an
algebra from ${TLb}_{7}.$
\begin{pro}\label{T2}\emph{(}Isomorphism criterion for ${TLb}_{7}$\emph{)}\\ Two filiform Leibniz algebras $L(C')$ and $L(C),$ where
$C'=(c'_{00},c'_{01},c'_{11},c'_{12},c'_{13},c'_{14},c'_{23})$ and
$C=(c_{00},c_{01},c_{11},c_{12},c_{13},c_{14},c_{23}),$ from ${TLb}_7$ are isomorphic if and only if there exist $A_0,A_1,B_1,B_2,B_3\in \mathbb{C}$ such that  $A_0B_1\neq 0$ and
the following hold:
\begin{align}
c_{00}^{\prime} &= \frac {{A^{2}_{{0}}}c_{00}+A_{{0}}A_{{1}}c_{01}+{A^{2}_{{
1}}}c_{11}}{{A^{5}_{{0}}}B_{{1}}}, \label{E1}\\[1mm]
c_{01}^\prime &= \frac
{A_{{0}}c_{01}+2\,A_{{1}}c_{11}}{{A^{5}_{{0}}}},\label{E2} \\ 
c_{11}^\prime &= \frac{B_1c_{11}}{A^{5}_0},\label{E3} \\
c_{12}^{\prime} &= \frac
{B_{{1}}c_{12}}{{A^{2}_{{0}}}},\label{E4} \\
c_{13}^{\prime} &= \frac{B_{1}\,c_{13}+B_{2}\,c_{12}}
{A^{3}_0},\label{E5} \\
c_{14}^{\prime} &= \frac{1}{A_{0}^{5}B_{1}}\left\{A_{0}B_{1}^{2}c_{14}+B_{2}(A_{0}B_{1}c_{13}+A_{1}B_{1}c_{12}c_{23}+A_{0}B_{2}c_{23}\right.\label{E6}\\&-\left.{} A_{1}B_{1}c_{12}^{2})-B_{3}(2A_{0}B_{1}c_{23}-A_{0}B_{1}c_{12})\right\}\nonumber \\
c_{23}^{\prime} &= \frac{B_{1}c_{23}}{A^{2}_{0}}\label{E7}
\end{align}
\begin{proof} 
 ``If'' part.\
 The equations in (\ref{E1}),\,(\ref{E2}),\,(\ref{E3}) occur due to {\bf Lemma \ref{L1}}.
Note that according  to {\bf Propositions \ref{p}} and {\bf \ref{p1}} the
adapted transformations in ${TLb}_{7}$ can be taken as follows:
\begin{equation}\label{BCH}
\left\{ \begin{array}{ll} 
e_0'& = f(e_0)=A_0\,e_0 + A_1\,e_1,
\\[1mm]
e_1'& = f(e_1)=B_1\,e_1+B_2\,e_2+B_3\,e_3  \\[1mm]
e_{i}'& = f(e_{i})=[f(e_{i-1}), f(e_0)],  \quad 1\leq i\leq 6,
\end{array}\right.
\end{equation}
where $ A_0\,B_1 \neq 0$ or more precisely:
\begin{eqnarray}
e_{0}^{\prime }&=&A_{0}e_{0}+A_{1}e_{1},\nonumber \\[1mm]
e_{1}^{\prime
}&=&B_{1}\,e_{1}+B_{2}e_{2}+B_{3}e_{3},\nonumber\\[1mm]
e_{2}^{\prime
}&=&A_{0}B_1\,e_{2}+A_{0}B_2e_{3}+\left(A_{{0}}B_{{3}}-A_{{1}}B_{{2}}c_{12}\right)e_4-A_{{1}}
\left( B_{{2}}c_{13}+B_{{3}}c_{12} \right) e_5 \nonumber\\[1mm]
&&+A_{{1}} \left( B_{{1}}c_{11}-B_{{2}}c_{14}-B_{{3}}c_{{
13}} \right)
e_6,\nonumber \\[1mm]
e_{3}^{\prime}&=&A_{0}^2B_1\,e_{3}+\left({A^{2}_{{0}}}B_{{2}}-A_{0}{A_{1}}B_1c_{12}\right)e_4+({A^{2}_{0}}B_{3}-2\,A_{0}A_{1}B_{2}c_{12}-A_{0}A_{1}B_1c_{13})e_5- \nonumber \\&& A_{1} \left(-A_{1}B_{{2}}{c^{2}_{{12}}}+A_{1}B_{2} c_{12}c_{23}+A_{0}B_{1}c_{14}+2\,A_{0}B_{2}c_{13}+2\,A_{0}B_{3}c_{12}-A_{0}B_{3}c_{23} \right)e_6,\nonumber \\[1mm]
e_{4}^{\prime}&=
&A_{0}^3B_1\,e_{4}+({A^{3}_{{0}}}B_{{2}}-2\,{A^{2}_{{0}}}A_{{1}}B_{{1}}c_{12})e_5+(A^{3}_
{0}B_{3}-2\,A_{1}A^{2}_{{0}}B_{1}c_{13}+A_{0}A^{2}_{1}B_{1}c^{2}_{12}\nonumber\\&& A_{0}A_{1}^{2}B_{1}c_{12}c_{23}-3\,A_{1}
A_{0}^{2}B_{2}c_{12}+A^{2}_{0}A_{1}B_{2}c_{23})e_6,\nonumber\\[1mm]
e_{5}^{\prime
}&=& A^4_0\,B_1\,e_5+({A^{4}_{{0}}}B_{{2}}-3\,A_{{1}}{A^{3}_{{0}}}B_{{1}}c_{12}+A_{{1}}{A^{3}_{{0}}}B_{{1}}c_
{{23}}
)e_6,\nonumber\\[1mm]
e_{6}^{\prime }&=& A^5_0\,B_1\,e_6.\nonumber
\end{eqnarray}
Consider
\begin{align} \label{bc.dim07}
[e_{1}^{\prime},e_{2}^{\prime}]=& A_{0}\,B_{1}^{2} \left(c_{12}e_{4}+c_{13}e_{5}+c_{14}e_{6}\right)+A_{0}B_{1}B_{2} \left(c_{12}e_{5}+c_{13}e_{6}\right)
+B_{1}(A_{2}B_{1}c_{12}-A_{1}B_{2}c_{12}+\nonumber\\&
A_{0}B_{3}) \left(c_{12}-c_{23}\right)e_{6}+A_{0}\,B_{2}^{2}c_{23}e_{6}-A_{0}B_{1}B_{3}c_{23}e_{6}\nonumber\\&=c^\prime_{12}e_{4}^{\prime}+c^\prime_{13}e_{5}^{\prime}+c_{14}^\prime e_{6}^{\prime} 
\end{align}
But
\begin{align*}
c'_{12}e'_{4}&=c'_{12}\big(A_{0}^3B_1\,e_{4}+({A^{3}_{{0}}}B_{{2}}-2\,{A^{2}_{{0}}}A_{{1}}B_{{1}}c_{12})e_5+(A^{3}_
{0}B_{3}-2\,A_{1}A^{2}_{{0}}B_{1}c_{13}+A_{0}A^{2}_{1}B_{1}c^{2}_{12}\\& -A_{0}A_{1}^{2}B_{1}c_{12}c_{23}-3\,A_{1}
A_{0}^{2}B_{2}c_{12}+A^{2}_{0}A_{1}B_{2}c_{23})e_6\big)&\\
c'_{13}e'_{5}&=c'_{13}\big(A^4_0\,B_1\,e_5+({A^{4}_{{0}}}B_{{2}}-3\,A_{{1}}{A^{3}_{{0}}}B_{{1}}c_{12}+A_{{1}}{A^{3}_{{0}}}B_{{1}}c_
{{2,3}})e_6\big)&\\
c'_{14}e'_{6}&=c'_{14}\big(A^5_0\,B_1\,e_6\big).
\end{align*}
On comparing the coefficients of $e_4, e_5$ and $e_6$ with the corresponding coefficients in (\ref{bc.dim07}) after rearrangement, we get the equalities (\ref{E4}), (\ref{E5}) and (\ref{E6}),
respectively. The last equality follows from \begin{align}
[f(e_2),f(e_3)]&= c^\prime_{23} f(e_6) \nonumber\\[2mm]
A_{0}^{3}B_{1}^{2}c_{23}e_{6}&=c'_{23}A_{0}
^{5}B_{1}e_{6}. \nonumber\\  
c_{23}^{\prime}&={\frac{B_1\,c_{23}}{{A^2_0}}}.\nonumber
 \end{align}
''Only if'' part. Let the equalities (\ref{E1})--(\ref{E7}) hold. Then the base
change above is adapted and $L(C)$  is transformed to $L(C').$

Indeed,
\begin{align*}
[e'_0,e'_0]&=\left[A_{0}e_{0}+A_{1}e_{1},\ A_{0}e_{0}+A_{1}e_{1}\right]\\[1mm]
&={A^{2}_{{0}}}[e_0,e_0]+A_{0}A_{1}[e_0,e_1]+A_{0}A_{1}[e_1,e_0]+A^2_{1}[e_1,e_1]\\[1mm]
&=\left({A^{2}_{{0}}}c_{00}+A_{{0}}A_{{1}}c_{01}+{A^{2}_{{1}}}c_{11}\right)e_{{6}}\\&
=c'_{00}{A^{5}_{{0}}}B_{{1}}e_6=c'_{00}e'_6.
\end{align*}
$c'_{00}$ is obtained as in (\ref{E1}) after dividing the last equation by the value of $e'_6$.
\begin{align*}
[e'_0,e'_1]&=\left[A_{0}e_{0}+A_{1}e_{1},\ B_{1}e_{1}+B_{2}e_{2}+B_3e_3\right]\\[1mm]&
= A_{0}B_{1} \left(-e_{2}+c_{01}e_{6} \right)-A_{0}B_{2}e_{3}-A_{0}B_{3}e_{4}+A_{1}B_{1}c_{11}e_{6}+\\& A_{1}B_{2}(c_{12}e_{4}+c_{13}e_{5}+c_{14}e_{6})+A_{1}B_{3}(c_{12}e_{5}+c_{13}e_{6})\\[1mm]
&=-(A_{0}B_1\,e_{2}+A_{0}B_2e_{3}+\left(A_{0}B_{3}- A_{1}B_{2}c_{12}\right)e_4-A_{1}
\left(B_{2}c_{13}+B_{3}c_{12} \right)e_5+\\& A_{1} \left( B_{1}c_{11}-B_{2}c_{14}-B_{3}c_{13} \right)
e_6)+B_{1} \left(A_{0}c_{01}+2\,A_{1}c_{11}
\right)e_{6}\\
 &= -e'_2+A_{0}^{5} B_{1}c'_{01}e_{6}= -e'_{2}+c'_{01}e'_{6}.
 \end{align*}
 $c'_{01}$ is obtained as in (\ref{E2}) after substituting the value of $e'_2$, simplify and divide resulting equation by the value of $e'_{6}$.
\begin{align*}
[e'_1,e'_1]&=\left[B_{1}e_{1}+B_{2}e_{2}+B_3e_3,B_{1}e_{1}+B_{2}e_{2}+B_3e_3\right]\\[1mm]
 &=\left[B_{{1}}e_1,B_{{1}}e_1\right]=B^2_1c_{11}e_6=A^5_0B_1c'_{11}e_6=c'_{11}e'_6.
 \end{align*}
 $c'_{11}$ is obtained as in (\ref{E3}) after dividing the last equation by the value of $e'_6$.\\ Other equations in {\bf Proposition \ref{T2}} 
are obtained similarly.
\end{proof}
 \end{pro}

\section{Isomorphism classes in $TLb_{7}$}

In this subsection, we give the list of isomorphism classes in
$TLb_{7}$. For simplification purpose, we introduce the function $\chi_{1}(X)=4\,x_{00}x_{11}-x_{01}^{2}$ where $X=(x_{00},x_{01},x_{11},x_{12},x_{13},x_{14},x_{23})$. By denoting $\Delta_{1}=\chi_{1}(C)\; \mbox{and}\; \Delta_{1}'=\chi_{1}(C')$, we present $TLb_{7}\ $ as a union of the following
subsets:\\
\indent$U_{7}^{1}=\{L(C)\in TLb_{7}\ :c_{23}\neq 0,\ c_{11} \neq 0,\ c_{12}\neq 0\};$

$U_{7}^{2}=\{L(C)\in TLb_{7}\ : c_{23}\neq 0,\ c_{11}\neq 0,\
c_{12}=0,\Delta_{1}\neq 0\};$

$U_{7}^{3}=\{L(C)\in TLb_{7}\ : c_{23}\neq 0,\
c_{11}\neq 0,c_{12}=\Delta_{1}=0\};$

$U_{7}^{4}=\{L(C)\in TLb_{7}\ : c_{23}\neq 0,\
c_{11}=0,c_{01}\neq 0,c_{12}\neq 0\};$

$U_{7}^{5}=\{L(C)\in TLb_{7}\ : c_{23}\neq 0,c_{11}=0,c_{01}\neq 0,c_{12}=0,c_{13}\neq 0\};$

$U_{7}^{6}=\{L(C)\in TLb_{7}\ : c_{23}\neq 0,\
c_{11}=0,c_{01}\neq 0,c_{12}=c_{13}=0\};$

$U_{7}^{7}=\{L(C)\in TLb_{7}\ : c_{23}\neq 0,\
c_{11}=c_{01}=0,c_{00}\neq 0,c_{12}\neq 0\};$

$U_{7}^{8}=\{L(C)\in TLb_{7}\ :c_{23}\neq 0,c_{11}=c_{01}=0,c_{00}\neq 0,c_{12}=0, c_{13}\neq 0\};$

$U_{7}^{9}=\{L(C)\in TLb_{7}\ :c_{23}\neq 0,c_{11}=c_{01}=0,c_{00}\neq 0,c_{12}=0,c_{13}=0\};$

$U_{7}^{10}=\{L(C)\in TLb_{7}\ :c_{23}\neq 0,c_{11}=c_{01}=c_{00}=0,c_{12}\neq 0,\};$

$U_{7}^{11}=\{L(C)\in TLb_{7}\ :c_{23}\neq 0,c_{11}=c_{01}=c_{00}=c_{12}=0,c_{13}\neq 0\};$

$U_{7}^{12}=\{L(C)\in TLb_{7}\ :c_{23}\neq 0, c_{11}=c_{01}=c_{00}=c_{12}=c_{13}=0\};$

$U_{7}^{13}=\{L(C)\in TLb_{7}\ :c_{23}=0,
c_{11}\neq 0,c_{12}\neq 0\};$

$U_{7}^{14}=\{L(C)\in TLb_{7}\ :c_{23}=0,
c_{11}\neq 0,c_{12}=0,\Delta_{1}\neq 0\};$

$U_{7}^{15}=\{L(C)\in TLb_{7}\ :c_{23}=0,
c_{11}\neq 0,c_{12}=\Delta_{1}=0,c_{13}\neq 0\};$

$U_{7}^{16}=\{L(C)\in TLb_{7}\ :c_{23}=0,c_{11}\neq 0,c_{12}=\Delta_{1}=c_{13}=0,c_{14}\neq 0\};$

$U_{7}^{17}=\{L(C)\in TLb_{7}\ :c_{23}=0,\
c_{11}\neq 0,\ c_{12}=\Delta_{1}=c_{13}=c_{14}=0\};$

$U_{7}^{18}=\{L(C)\in TLb_{7}\ :c_{23}=c_{11}=0,c_{01}\neq 0,c_{12}\neq 0\};$

$U_{7}^{19}=\{L(C)\in TLb_{7}\ :c_{23}=c_{11}=0,c_{01}\neq 0,c_{12}=0,c_{13}\neq 0\};$

$U_{7}^{20}=\{L(C)\in TLb_{7}\
:c_{23}=c_{11}=0,c_{01}\neq 0,c_{12}=c_{13}=0,c_{14}\neq 0\};$

$U_{7}^{21}=\{L(C)\in TLb_{7}\
:c_{23}=c_{11}=0,c_{01}\neq 0,c_{12}=c_{13}=c_{14}=0\};$

$U_{7}^{22}=\{L(C)\in TLb_{7}\
:c_{23}=c_{11}=c_{01}=0,\ c_{00}\neq 0,\
c_{12}\neq 0\};$

$U_{7}^{23}=\{L(C)\in TLb_{7}\
:c_{23}=c_{11}=c_{01}=0,\ c_{00}\neq 0,\
c_{12}=0,c_{13}\neq 0\};$

$U_{7}^{24}=\{L(C)\in TLb_{7}\
:c_{23}=c_{11}=c_{01}=0,\ c_{00}\neq 0,\
c_{12}=0,c_{13}=0,c_{14}\neq 0\};$

$U_{7}^{25}=\{L(C)\in TLb_{7}\
:c_{23}=c_{11}=c_{01}=0,c_{00}\neq 0,
c_{12}=c_{13}=c_{14}=0\};$

$U_{7}^{26}=\{L(C)\in TLb_{7}\
:c_{23}=c_{11}=c_{01}=c_{00}=0,\
c_{12}\neq 0\}.$

$U_{7}^{27}=\{L(C)\in TLb_{7}\
:c_{23}=c_{11}=c_{01}=c_{00}=c_{12}=0,c_{13}\neq 0,c_{14}\neq 0\}.$

$U_{7}^{28}=\{L(C)\in TLb_{7}\
:c_{23}=c_{11}=c_{01}=c_{00}=c_{12}=c_{13}\neq 0,c_{14}=0\}.$

$U_{7}^{29}=\{L(C)\in TLb_{7}\
:c_{23}=c_{11}=c_{01}=c_{00}=c_{12}=c_{13}=0,c_{14}\neq 0\}.$

$U_{7}^{30}=\{L(C)\in TLb_{7}\
:c_{23}=c_{11}=c_{01}=c_{00}=c_{12}=c_{13}=c_{14}=0\}.$\\
The adapted group action on $TLb_7$ is induced on each of these disjoint subsets. The next proposition give isomorphism criterion for one subset from each type, that is: from infinitely many orbits case and from single orbits case, since other cases can be proved similarly.
\begin{pro}\label{p2} \emph{}
\begin{enumerate}
\item[\em{1.}]
\begin{enumerate}
\item[\em{(a)}]  Two algebras $L(C')$ and $L(C)$ from $U_{7}^{1}$ are
isomorphic, if and only if \\ $$ \left(\frac
{c'_{23}}{c'_{11}}\right)^{8} \Delta_{1}'^{3}=\left(\frac
{c_{23}}{c_{11}}\right)^{8} \Delta_{1}^{3} \\[2mm]
$$ \mbox{and} $$ \frac {c^\prime_{12}}{c'_{23}}={\frac {c_{12}}{c_{23}}} \ $$
 \item[\em{(b)}] For any $\lambda_1,  \lambda_2, \in \mathbb{C}, $  there exists $L(C)\in
U_{7}^{1}$ such that 
$$ \left(\frac
{c_{23}}{c_{11}}\right)^{8}\,\Delta_{1}^{3}=\lambda_1,
\ \ {\frac {c_{12}}{c_{23}}}=\lambda_2.$$ {\em Then orbits in
$U_{7}^{1}$ can be parametrized as} $L\left(\lambda_1, 0, 1,
\lambda_2,0,0,1\right), \ \lambda_1,
\lambda_2 \in \mathbb{C}.$
\end{enumerate}
\item[\em{2}] The subset $ U_{7}^{6}$ is a single orbit with
the representative \ $ L(0,1,0,0,0,0,1).$
\end{enumerate}
\begin{proof}
\begin{enumerate}
\item
\begin{enumerate}
\item The ``If'' part of the proposition is due to {\bf Proposition \ref{T2}} if one substitutes the expressions for $c_{00}',$ $c_{01}',$ $c_{11}',$ $c_{12}',$
$c_{13}',$ $c_{23}'$ into
$$\left(\frac
{c_{23}'}{c_{11}'}\right)^{8}\ \Delta_{1}'^{3},\ \
\mbox{and} \ \ \frac {c^\prime_{12}}{c'_{23}},\ \ \mbox{respectively.}$$

`` Only if '' part. \ Let the equalities  $$  \left(\frac
{c'_{23}}{c'_{11}}\right)^{8}\ \Delta_{1}'^{3}=\left(\frac
{c_{23}}{c_{11}}\right)^{8}\ \Delta_{1}^{3}$$
and $$ \frac {c^\prime_{12}}{c'_{23}}={\frac {c_{12}}{c_{23}}}$$ hold.
Consider the base change (\ref{BCH}) in the proof of {\bf Proposition \ref{T2}} with
$$ A_1={\frac {-A_{{0}}c_{01}}{2\,c_{11}}}, \
B_1={\frac {{A_{{0}}}^{2}}{c_{{23}}}},$$ and 
\begin{align*}
 B_3=\frac{A_0B_{1}^{2}c_{14}+B_{2}(A_0B_1c_{13} + A_1B_1c_{12}c_{23}+ A_0B_2c_{23}-A_1B_1c_{12}^{2})}{2A_0B_1c_{23}-A_0B_1c_{12}}
\end{align*}
This base change transforms $L(C)$ into
$$ L\left(\left(\frac
{c_{23}}{c_{11}}\right)^{8} \Delta_{1}^{3},0,1,{\frac
{c_{12}}{c_{23}}},0,0,1\right).$$ An analogous
base change with ``prime'' transforms $L(C')$ into $$
L\left(\left(\frac
{c'_{23}}{c'_{11}}\right)^{8}\Delta_{1}'^{3},0,1,{\frac
{{c'}_{{12}}}{{c'}_{23}}},0,0,1\right).$$ Since
$$ \left(\frac
{c'_{23}}{c'_{11}}\right)^{8}\Delta_{1}'^{3}=\left(\frac
{c_{23}}{c_{11}}\right)^{8}\Delta_{1}^{3} $$ \ \ and
$$ \frac {c^\prime_{12}}{c'_{23}}={\frac {c_{12}}{c_{23}}},$$ the algebras are
isomorphic. \qquad Obvious.
 \end{enumerate}
\item To prove it, it is enough to put appropriate values of $ A_{0},A_{1},B_{1}$ in the base change (\ref{BCH}).\\
Then, it is easy to see that $A^{4}_{0}=c_{01}, A_{1}=-A_{0}\frac{c_{00}}{c_{01}}, B^{2}_{1}=\frac{c_{01}}{c^{2}_{23}}$ and $ B_{2}, B_{3}\in \mathbb{C}.$
 \end{enumerate}
\end{proof}
\end{pro}
\section{Isomorphism criterion for ${TLb}_{8}$}
As noted in the introductory section and corroborated in Section \ref{sec:2}, an algebra $L$ from ${TLb}_{8}$ is represented as a one dimensional Leibniz central extension\,($C(L)=<e_7>$)
 of the filiform Lie algebra with law $\mu_7+\Psi,$ where $\Psi$ is the appropriate cocycle on the adapted basis $\{e_0,e_1,...,e_7\},$ the class ${TLb}_{8}$ is represented as follows:
 \begin{center}
 $TLb_{8}=\left\{
 \begin{array}{lll}
 [e_{i},e_{0}]=e_{i+1}, \qquad \qquad \qquad \qquad \qquad \quad \qquad \qquad  1\leq i\leq
 {6},   \\[0mm]
 [e_{0},e_{i}]=-e_{i+1}, \qquad \qquad \qquad \qquad \qquad \qquad \qquad \ 2\leq i\leq {6},   \\[0mm]
 [ e_{0},e_{0}]=b_{00}e_{7},  \\[0mm]
 [ e_{0},e_{1}]=-e_{2}+b_{01}e_{7}, \\[0mm]
 [ e_{1},e_{1}]=b _{11}e_{7}, \\[0mm]
 [ e_{1},e_{2}]=-[e_2,e_1]=a_{14}e_4+a_{15}e_5+a_{16}e_{6}+b_{12}e_{7} , \\[0mm]
 [e_{1},e_{3}]=-[e_{3},e_{1}]=a_{14}\,e_5+a_{15}e_{6}+b_{13}e_{7}, \\[0mm]
  [e_{1},e_{4}]=-[e_{4},e_{1}]=(a_{14}-a_{26})e_6+b_{14}e_{7}, \\[0mm]
  [e_{1},e_{5}]=-[e_{5},e_{1}]=b_{15}e_{7},\\[0mm]
  [e_{2},e_{3}]=-[e_{3},e_{2}]=a_{26}e_6+b_{23}e_{7},\\[0mm]
  [e_{2},e_{4}]=-[e_{4},e_{2}]=b_{24}e_{6}.\\[0mm]
  [e_{i},e_{7-i}]=-[e_{7-i},e_{i}]=(-1)^{i}b_{34}e_{7} \qquad \qquad \qquad \;  1\leq i\leq {6}
 \end{array}%
 \right. $ 
 \end{center}
The interrelations between the structure constants $a_{ij}$ and $b_{st}$ in the table above are listed as follows.
\begin{lem}
 The structure constants of algebras from ${TLb}_{8}$ satisfy the following constraints: 
 \begin{enumerate}
 \item[{\em 1.}] $b_{13}=a_{16}$, 
 \item[{\em 2.}] $b_{14}=a_{15}-b_{23}$
 \item[{\em 3.}] $b_{24}=a_{26}$
 \item[{\em 4.}] $b_{15}=a_{14}-2a_{26}$
 \item[{\em 5.}] $b_{34}(a_{26}+2a_{14})=0$
\end{enumerate}
\begin{proof}
First and  second  relations occur if we apply the Leibniz identity to
the triples $\{e_0,e_1,e_2\},\
\{e_0,e_1,e_3\}$ of the basis vectors  respectively, and
the third, fourth and fifth relations are obtained as a result of application of the Leibniz identity to $\{e_0,e_2,e_3\},\{e_0,e_1,e_4\},$ and $\{e_1,e_2,e_3\}$ respectively. 
\end{proof}
\end{lem}
Further unifying the table of multiplication above we rewrite it via new notations
$c_{00},c_{01},c_{11},c_{12},c_{13},c_{14},c_{15},c_{23},c_{24},c_{34}$ for the structure constants:
 \begin{center}
 $TLb_{8}=\left\{
 \begin{array}{lll}
 [e_{i},e_{0}]=e_{i+1}, \qquad \qquad \qquad\qquad \qquad \qquad \qquad \quad  1\leq i\leq
 {6},   \\[0mm]
 [e_{0},e_{i}]=-e_{i+1}, \qquad \qquad \qquad\qquad \qquad \qquad \qquad \ 2\leq i\leq {6},   \\[0mm]
 [ e_{0},e_{0}]=c_{00}e_{7},  \\[0mm]
 [ e_{0},e_{1}]=-e_{2}+c_{01}e_{7}, \\[0mm]
 [ e_{1},e_{1}]=c_{11}e_{7}, \\[0mm]
 [e_{1},e_{2}]=-[e_2,e_1]=c_{12}e_4+c_{13}e_5+c_{14}e_{6}+c_{15}e_{7} , \\[0mm]
 [e_{1},e_{3}]=-[e_{3},e_{1}]=c_{12}e_5+c_{13}e_{6}+c_{14}e_{7}, \\[0mm]
  [e_{1},e_{4}]=-[e_{4},e_{1}]=(c_{12}-c_{23})e_6+(c_{13}-c_{24})e_{7}, \\[0mm]
  [e_{1},e_{5}]=-[e_{5},e_{1}]=(c_{12}-2c_{23})e_{7},\\[0mm]
  [e_{2},e_{3}]=-[e_{3},e_{2}]=c_{23}e_6+c_{24}e_{7},\\[0mm]
  [e_{2},e_{4}]=-[e_{4},e_{2}]=c_{23}e_{7}.\\[0mm]
  [e_{i},e_{7-i}]=-[e_{7-i},e_{i}]=(-1)^{i}c_{34}e_{7} \qquad \qquad \qquad \;  1\leq i\leq {6}
 \end{array}%
 \right. $ 
 \end{center}
 An algebra from ${TLb}_{8}$ with the structure constants
 $c_{00},c_{01},c_{11},c_{12},c_{13},c_{14},c_{15},c_{23},c_{24},c_{34}$ is
 denoted by $L(C),$ where
 $C=(c_{00},c_{01},c_{11},c_{12},c_{13},c_{14},c_{15},c_{23},c_{24},c_{34}).$
 The next theorem represents the action of the adapted base change on the structure constants
$c_{00},c_{01},c_{11},c_{12},c_{13},c_{14},c_{15},c_{23},c_{24},c_{34}$ of algebra from ${TLb}_{8}.$
 \begin{pro}\label{T3}\emph{(}Isomorphism criterion for ${TLb}_{8}$\emph{)}\\ Two filiform Leibniz algebras $L(C')$ and $L(C),$ where
 $C'=(c'_{00},c'_{01},c'_{11},c'_{12},c'_{13},c'_{14},c'_{23})$ and
 $C=(c_{00},c_{01},c_{11},c_{12},c_{13},c_{14},c_{23}),$ from ${TLb}_8$ are isomorphic if and only if there exist $A_0,A_1,B_1,B_2,B_3\in \mathbb{C}$ such that  $A_0B_1\neq 0$ and the following hold:
\begin{align} 
 c_{00}^\prime&=\frac {{A^{2}_{{0}}}c_{00}+A_{{0}}A_{{1}}c_{01}+{A^{2}_{{
 1}}}c_{11}}{{A^{5}_{{0}}}B_{{1}}(A_{0}+A_{1}c_{34})}, \label{D2}\\[1mm] 
 c_{01}^\prime&=\frac
 {A_{{0}}c_{01}+2\,A_{{1}}c_{11}}{{A^{5}_{{0}}}(A_{0}+A_{1}c_{34})}, \label{D3}\\ 
 c_{11}^\prime&=\frac{B_1c_{11}}{A^{5}_{0}(A_{0}+A_{1}c_{34})}, \label{D4}\\
 c_{12}^{\prime}&=\frac
 {B_{{1}}c_{12}}{{A^{2}_{{0}}}}, \label{D5}\\
 c_{13}^{\prime}&=\frac{B_{1}\,c_{13}+B_{2}\,c_{12}} {A^{3}_0}, \label{D6}\\
 c_{14}^{\prime}&=\frac{1}{A_{0}^{5}B_{1}}\left\{A_{0}B_{1}^{2}c_{14}+B_{2}(A_{0}B_{1}c_{13}+A_{1}B_{1}c_{12}c_{23}+A_{0}B_{2}c_{23}-A_{1}B_{1}c_{12}^{2})\label{D7}\right.\\&-\left. {}B_{3}(2A_{0}B_{1}c_{23}-A_{0}B_{1}c_{12})\right\}\notag \\
 c_{15}^{\prime}&=\frac{1}{A_{0}^{5}B_{1}(A_{0}+A_{1}c_{34})}\left\{A_{0}B_{1}^{2}c_{15}+(A_{0}B_{1}B_{2}+A_{1}B_{1}B_{2})c_{14}+(A_{0}B_{1}B_{3}-2A_{1}B_{1}B_{2}c_{12}\notag\right.\\&+\left.  {}2A_{1}B_{1}B_{2}c_{23}+A_{1}B_{1}B_{3}c_{34})c_{13}+(A_{0}B_{1}B_{4}+A_{1}B_{1}B_{2}c_{24}-A_{1}B_{1}B_{3}c_{12}+2A_{1}B_{1}B_{3}c_{23}\notag\right.\\&+\left. {}A_{1}B_{1}B_{4}c_{34}-A_{1}B^{2}_{2}c_{23})c_{12}+(A_{0}B_{2}B_{3}-3A_{0}B_{1}B_{4}-A_{1}B_{1}B_{4}c_{34})c_{23}+(A_{0}B^{2}_{2}\notag\right.\\&-\left. {}2A_{0}B_{1}B_{3})c_{24}+(2A_{0}B_{2}B_{4}-A_{0}B^{2}_{3}-2A_{0}B_{1}B_{5})c_{34}\right\} \label{D8}\\
 c_{23}^{\prime}&=\frac{B_{1}c_{23}}{A^{2}_{0}}\label{D9}\\
 c_{24}^{\prime}&=\frac{1}{A^{5}_{0}(A_{0}+A_{1}c_{34})}\left\{A^{3}_{0}B^{2}_{1}c_{24}+(A^{3}_{0}B_{1}B_{2}-A^{2}_{0}A_{1}B^{2}_{1}c_{12})c_{23}+(2A^{3}_{0}B_{1}B_{3}-A^{3}_{0}B^{2}_{2}\notag\right.\\&-\left.2A^{2}_{0}A_{1}B_{1}B_{2}c_{12}-A^{2}_{0}A_{1}B^{2}_{1}c_{13})c_{34}\right\} \label{D10}\\
 c_{34}^{\prime}&=\frac{B_{1}c_{34}}{A_{0}+A_{1}c_{34}}\label{D11}
 \end{align}
 \newpage
 \begin{proof} ``If'' part.\
 The equations (\ref{D2}),\,(\ref{D3}),\,(\ref{D4}) occur due to {\bf Lemma \ref{L1}}.
Note that according  to {\bf Propositions \ref{p}} and {\bf \ref{p1}}, the adapted transformations in ${TLb}_{8}$ can be taken as follows: 
\begin{equation}
\left\{ \begin{aligned} \label{BC}
e_0'& = f(e_0)=A_0\,e_0 + A_1\,e_1,
\\[1mm]
e_1'& = f(e_1)=B_1\,e_1+B_2\,e_2+B_3\,e_3+B_4\,e_4+B_5\,e_5,  \\[1mm]
e_{i}'& = f(e_{i})=[f(e_{i-1}), f(e_0)],  \quad 1\leq i\leq 7,
\end{aligned}\right.
\end{equation}
where $ A_0\,B_1 \neq 0,$ that is, \vspace{-1cm}
\begin{align}
e_{2}^{\prime
}&=A_{0}B_1e_{2}+A_{0}B_2e_{3}+\left(A_{{0}}B_{{3}}-A_{{1}}B_{{2}}c_{12}\right)e_4-\left(A_{0}B_{4}+A_{1}B_{2}c_{13}+A_{1}B_{3}c_{12}\right)e_5 \nonumber\\[1mm]
&+(A_{0}B_{5}-A_{1}B_{2}c_{14}-A_{1}B_{3}c_{13}+A_{1}B_{4}c_{23}-A_{1}B_{4}c_{12})e_{6}+\left(A_{1}B_{1}c_{11}\nonumber\right.\\&-\left. {}A_{1}B_{2}c_{15}-A_{1}B_{3}c_{14}+A_{1}B_{4}c_{24}-A_{1}B_{4}c_{13}+2A_{1}B_{5}c_{23}-A_{1}B_{5}c_{12}\right)e_7,\nonumber \\[1mm]
e_{3}^{\prime}&=A_{0}^2B_1e_{3}+\left({A^{2}_{0}}B_{2}-A_{0}A_{1}B_1c_{12}\right)e_4+(A^{2}_{0}B_{3}-2A_{0}A_{1}B_{2}c_{12}-A_{0}A_{1}B_1c_{13})e_5\nonumber \\&+ \left(A^{2}_{0}B_{4}-2A_{0}A_{1}B_{2}c_{13}-2A_{0}A_{1}B_{3}c_{12}-A_{0}A_{1}B_{1}c_{14}+A_{0}A_{1}B_{3}c_{23}\nonumber\right.\\&-\left. {}A^{2}_{1}B_{2}c_{12}c_{23}+A^{2}_{1}B_{2}c^{2}_{12}\right)e_6+(A^{2}_{0}B_{5}-2A_{0}A_{1}B_{2}c_{14}-A_{0}A_{1}B_{3}c_{13}+3A_{0}A_{1}B_{4}c_{23}\nonumber\\&-2A_{0}A_{1}B_{4}c_{12}-A_{0}A_{1}B_{1}c_{15}+A_{0}A_{1}B_{3}c_{24}+2A^{2}_{1}B_{2}c_{12}c_{13}-2A^{2}_{1}B_{2}c_{13}c_{23}-2A^{2}_{1}B_{3}c_{12}c_{23}\nonumber\\&+A^{2}_{1}B_{3}c^{2}_{12}+A_{0}A_{1}B_{5}c_{34}-A^{2}_{1}B_{2}c_{14}c_{34}-A^{2}_{1}B_{3}c_{12}c_{34}+A^{2}_{1}B_{4}c_{23}c_{34}-A^{2}_{1}B_{4}c_{12}c_{34})e_{7},\nonumber \\[1mm]
e_{4}^{\prime}&=A_{0}^3B_{1}e_{4}+({A^{3}_{0}}B_{{2}}-2A^{2}_{0}A_{1}B_{1}c_{12})e_5+(A^{3}_{0}B_{3}-2A^{2}_{0}A_{1}B_{2}c_{12}-2A^{2}_{0}A_{1}B_{1}c_{13})e_{6}\nonumber\\&+(A^{3}_{0}B_{4}-2A^{2}_{0}A_{1}B_{2}c_{13}-2A^{2}_{0}A_{1}B_{3}c_{12}-2A^{2}_{0}A^{2}_{1}B_{1}c_{14}+A^{2}_{0}A_{1}B_{3}c_{23}-A_{0}A^{2}_{1}B_{2}c_{12}c_{23}\nonumber\\&+A_{0}A^{2}_{1}B_{2}c^{2}_{12}+2A^{3}_{0}A_{1}B_{3}c_{23}-A^{3}_{0}A_{1}B_{2}c_{12}-2A^{2}_{0}A^{2}_{1}B_{1}c_{12}c_{23}+A^{2}_{0}A^{2}_{1}B_{1}c^{2}_{12}+A^{3}_{0}A_{1}B_{3}c_{34}\nonumber\\&-A^{2}_{0}A^{2}_{1}B_{1}c_{13}c_{34}-2A^{2}_{0}A^{2}_{1}B_{2}c_{12}c_{34})e_{7},
\nonumber \\[1mm]
e_{5}^{\prime}&=A^{4}_{0}B_{1}e_5+(A^{4}_{0}B_{2}+A^{3}_{0}A_{1}B_{1}c_{23}-3A^{3}_{0}A_{1}B_{1}c_{12})e_{6}+(A^{4}_{0}B_{3}-3A^{3}_{0}A_{1}B_{1}c_{12}-2A^{3}_{0}A_{1}B_{1}c_{13}\nonumber\\&+A^{3}_{0}A_{1}B_{1}c_{24}+2A^{3}_{0}A_{1}B_{2}c_{23}-4A^{2}_{0}A^{2}_{1}B_{1}c_{12}c_{23}-2A^{2}_{0}A^{2}_{1}B_{1}c^{2}_{12}+A^{3}_{0}A_{1}B_{3}c_{34}\nonumber\\&-2A^{2}_{0}A^{2}_{1}B_{2}c_{12}c_{34}-2A^{2}_{0}A^{2}_{1}B_{1}c_{13}c_{34})e_7,\nonumber\\[1mm]
e_{6}^{\prime}&=A^{5}_{0}B_{1}e_{6}+(A^{5}_{0}B_{2}+3A^{4}_{0}A_{1}B_{1}c_{23}-4A^{4}_{0}A_{1}B_{1}c_{12}+A^{4}_{0}A_{1}B_{2}c_{34}+A^{3}_{0}A^{2}_{1}B_{1}c_{23}c_{34}\nonumber\\&-3A^{3}_{0}A^{2}_{1}B_{1}c_{13}c_{34})e_{7}\nonumber\\&
e_{7}^{\prime}=A^{5}_{0}B_{1}(A_{0}+A_{1}c_{34})e_{7}.\nonumber
\end{align}
Consider,
\begin{align} \label{bc.dim08}
[e_{1}^{\prime},e_{2}^{\prime}]=& A_{0}B_{1}^{2}c_{12}e_{4}+\left(A_{0}B^{2}_{1}c_{13}+A_{0}B_{1}B_{2}c_{12}\right)e_{5}+\left(A_{0}B^{2}_{1}c_{14}+A_{0}B_{1}B_{3}c_{13}\nonumber\right.\\&+\left. {}(A_{0}B_{1}B_{3}+A_{1}B_{1}B_{2}c_{23})c_{12}+(A_{0}B^{2}_{2}-2A_{0}B_{1}B_{3})c_{23}-A_{1}B_{1}B_{2}c^{2}_{12}\right)e_{6}\nonumber\\&+\left(A_{0}B^{2}_{1}c_{15}+(A_{0}B_{1}B_{2}+A_{1}B_{1}B_{2}c_{34})c_{14}+(A_{0}B_{1}B_{3}-2A_{1}B_{1}B_{2}c_{12}+2A_{1}B_{1}B_{2}c_{23}\nonumber\right.\\&+\left. {}A_{1}B_{1}B_{3}c_{34})c_{13}+(A_{0}B_{1}B_{4}+A_{1}B_{1}B_{2}c_{24}-A_{1}B_{1}B_{3}c_{12}+2A_{1}B_{1}B_{3}c_{23}\nonumber\right.\\&+\left. {}A_{1}B_{1}B_{4}c_{34}-A_{1}B^{2}_{2}c_{23}-A_{1}B_{2}B_{3}c_{34}+A_{1}B_{2}B_{3}c_{34})c_{12}+(A_{0}B_{2}B_{3}-3A_{0}B_{1}B_{4}\nonumber\right.\\&-\left. {}A_{1}B_{1}B_{4}c_{34})c_{23}+(A_{0}B^{2}_{2}-2A_{0}B_{1}B_{3})c_{24}+(2A_{0}B_{2}B_{4}-A_{0}B^{2}_{3}-2A_{0}B_{1}B_{5}c_{34}) \right)e_{7}\nonumber\\&=c^\prime_{12}e_{4}^{\prime}+c^\prime_{13}e_{5}^{\prime}+c_{14}^\prime e_{6}^{\prime}+c_{15}^\prime e_{7}^{\prime} 
\end{align}

But,
\begin{align*}
c'_{12}e'_{4}&=c'_{12}\Big(A_{0}^3B_{1}e_{4}+({A^{3}_{0}}B_{{2}}-2A^{2}_{0}A_{1}B_{1}c_{12})e_5+(A^{3}_{0}B_{3}-2A^{2}_{0}A_{1}B_{2}c_{12}-2A^{2}_{0}A_{1}B_{1}c_{13})e_{6}\nonumber\\&+(A^{3}_{0}B_{4}-2A^{2}_{0}A_{1}B_{2}c_{13}-2A^{2}_{0}A_{1}B_{3}c_{12}-2A^{2}_{0}A^{2}_{1}B_{1}c_{14}+A^{2}_{0}A_{1}B_{3}c_{23}-A_{0}A^{2}_{1}B_{2}c_{12}c_{23}\nonumber\\&+A_{0}A^{2}_{1}B_{2}c^{2}_{12}+2A^{3}_{0}A_{1}B_{3}c_{23}-A^{3}_{0}A_{1}B_{2}c_{12}-2A^{2}_{0}A^{2}_{1}B_{1}c_{12}c_{23}+A^{2}_{0}A^{2}_{1}B_{1}c^{2}_{12}+A^{3}_{0}A_{1}B_{3}c_{34}\nonumber\\&-A^{2}_{0}A^{2}_{1}B_{2}c_{12}c_{34}-A^{2}_{0}A^{2}_{1}B_{1}c_{13}c_{34}-A^{2}_{0}A_{1}B_{2}c_{12}c_{34})e_{7}\Big),
\nonumber&\\
c'_{13}e'_{5}&=c'_{13}\Big(A^{4}_{0}B_{1}e_5+(A^{4}_{0}B_{2}+A^{3}_{0}A_{1}B_{1}c_{23}-3A^{3}_{0}A_{1}B_{1}c_{12})e_{6}+(A^{4}_{0}B_{3}-3A^{3}_{0}A_{1}B_{1}c_{12}\nonumber\\&-2A^{3}_{0}A_{1}B_{1}c_{13}+A^{3}_{0}A_{1}B_{1}c_{24}+2A^{3}_{0}A_{1}B_{2}c_{23}-4A^{2}_{0}A^{2}_{1}B_{1}c_{12}c_{23}-2A^{2}_{0}A^{2}_{1}B_{1}c^{2}_{12}\nonumber\\&+A^{3}_{0}A_{1}B_{3}c_{34}-2A^{2}_{0}A^{2}_{1}B_{2}c_{12}c_{34}-2A^{2}_{0}A^{2}_{1}B_{1})e_7\Big),\nonumber&\\
c'_{14}e'_{6}&=c'_{14}\Big(A^{5}_{0}B_{1}e_{6}+(A^{5}_{0}B_{2}+3A^{4}_{0}A_{1}B_{1}c_{23}-4A^{4}_{0}A_{1}B_{1}c_{12}+A^{4}_{0}A_{1}B_{2}c_{34}+A^{3}_{0}A^{2}_{1}B_{1}c_{23}c_{34}\nonumber\\&-3A^{3}_{0}A^{2}_{1}B_{1}c_{13}c_{34})e_{7}\Big),\nonumber&\\
c'_{15}e'_{7}&=c'_{15}\Big(A^{5}_{0}B_{1}(A_{0}+A_{1}c_{34})e_{7}\Big)
\end{align*}
On comparing the coefficients of $e_4, e_5, e_6$ and  $e_7$ with the corresponding coefficients in (\ref{bc.dim08}), we obtain the equalities (\ref{D5}), (\ref{D6}), (\ref{D7}) and (\ref{D8}),
respectively. The equalities (\ref{D9}) and (\ref{D10}) follows from
\begin{align}\label{89}
[f(e_1),f(e_4)]&= c^\prime_{23} e^\prime_{6}+c^\prime_{24}e^\prime_{7} \nonumber\\[2mm]
[f(e_1),f(e_4)]&=A_{0}B^{2}_{1}c_{23}e_{6}+\Big(A^{3}_{0}B^{2}_{1}c_{24}+(A^{3}_{0}B_{1}B_{2}-2A^{2}_{0}A_{1}B^{2}_{1}c_{12})c_{23}+A^{2}_{0}A_{1}B^{2}_{1}c^{2}_{12}\nonumber\\&-(2A^{3}_{0}B_{1}B_{3}-2A^{2}_{0}A_{1}B^{2}_{1}c_{12}-2A^{2}_{0}A_{1}B^{2}_{1}c_{13}+A^{2}_{0}A_{1}B_{1}B_{2}c_{12}-A^{3}_{0}B^{2}_{2})c_{34}\Big)e_{7}
\end{align}
\mbox{Since} 
\begin{align*}
c^\prime_{23} e^\prime_{6}&=c^\prime_{23}\Big(A^{5}_{0}B_{1}e_{6}+(A^{5}_{0}B_{2}+3A^{4}_{0}A_{1}B_{1}c_{23}-4A^{4}_{0}A_{1}B_{1}c_{12}+A^{4}_{0}A_{1}B_{2}c_{34}+A^{3}_{0}A^{2}_{1}B_{1}c_{23}c_{34}\nonumber\\&-3A^{3}_{0}A^{2}_{1}B_{1}c_{13}c_{34})e_{7}\Big) \ \mbox{and} \nonumber&\\
c^\prime_{24}e^\prime_{7}&=\Big(A^{5}_{0}B_{1}(A_{0}+A_{1}c_{34})e_{7}\Big).\nonumber
\end{align*}
We obtain the equalities (\ref{D9}) and (\ref{D10}) by comparing the coefficients of $e_{6}$ and $e_{7}$ with the corresponding coefficients in (\ref{89}) respectively. \\
The last equation follows by considering the bracket
\begin{align*}
[f(e_3),f(e_4)]&=A^{5}_{0}B^{2}_{1}c_{34}e_{7}=c'_{34}e'_{7}=c'_{34}(A^{5}_{0}B_{1}(A_{0}+A_{1}c_{34})e_{7}).
\end{align*}
``Only if'' part. Let the equalities (\ref{D2})--(\ref{D11}) hold. Then the base
change (\ref{BC}) is adapted and $L(C)$ is transformed to $L(C').$

Indeed,
\begin{align*}
[e'_0,e'_0]&=\left[A_{0}e_{0}+A_{1}e_{1},\ A_{0}e_{0}+A_{1}e_{1}\right]\\[1mm]
&={A^{2}_{{0}}}[e_0,e_0]+A_{0}A_{1}[e_0,e_1]+A_{0}A_{1}[e_1,e_0]+A^2_{1}[e_1,e_1]\\[1mm]
&=\left({A^{2}_{{0}}}c_{00}+A_{{0}}A_{{1}}c_{01}+{A^{2}_{{1}}}c_{11}\right)e_{{6}}\\&=c'_{00}{A^{5}_{{0}}}B_{{1}}e_6=c'_{00}e'_6.&\\
\end{align*}
\begin{align*}
[e'_0,e'_1]&=\left[A_{0}e_{0}+A_{1}e_{1},\ B_{1}e_{1}+B_{2}e_{2}+B_3e_3\right]\\[1mm]&
= A_{0}B_{1} \left(-e_{2}+c_{01}e_{6} \right)-A_{0}B_{2}e_{3}-A_{0}B_{3}e_{4}+A_{1}B_{1}c_{11}e_{6}+\\& A_{1}B_{2}(c_{12}e_{4}+c_{13}e_{5}+c_{14}e_{6})+A_{1}B_{3}(c_{12}e_{5}+c_{13}e_{6})\\[1mm]
&=-(A_{0}B_1\,e_{2}+A_{0}B_2e_{3}+\left(A_{0}B_{3}- A_{1}B_{2}c_{12}\right)e_4-A_{1}
\left(B_{2}c_{13}+B_{3}c_{12} \right)e_5+\\& A_{1} \left( B_{1}c_{11}-B_{2}c_{14}-B_{3}c_{13} \right)
e_6)+B_{1} \left(A_{0}c_{01}+2\,A_{1}c_{11}
\right)e_{6}\\
 &= -e'_2+A_{0}^{5} B_{1}c'_{01}e_{6}= -e'_{2}+c'_{01}e'_{6}. &\\
[e'_1,e'_1]&=\left[B_{1}e_{1}+B_{2}e_{2}+B_3e_3,B_{1}e_{1}+B_{2}e_{2}+B_3e_3\right]\\[1mm]
 &=\left[B_{{1}}e_1,B_{{1}}e_1\right]=B^2_1c_{11}e_6=A^5_0B_1c'_{11}e_6=c'_{11}e'_6. 
\end{align*}
$c'_{01}$ and $c'_{11}$ are obtained as in (\ref{D2}) and (\ref{D3}) after substituting the value of $e'_2$ and $e'_6$, simplify and divide resulting equation by the value of $e'_{6}$ respectively.
The other parameters on the left-hand side of the equations in {\bf Proposition \ref{T3}} can be obtained similarly by considering $[e_{1}',e_{2}']$ and $[e_{1}',e_{4}']$.
\end{proof}
 \end{pro}
\section{Isomorphism classes in $TLb_{8}$}
In this subsection we give the list of isomorphism classes in
$TLb_{8}$. For simplification purpose, we introduce the following notations $\chi_{2}(X)=2\,x_{11}-x_{01}x_{34},$ $\chi_{3}(X)=x_{01}-x_{00}x_{34},$ and $\chi_{4}(X)=x_{12}x_{23}+x_{13}x_{34}$ where $X=(x_{00},x_{01},x_{11},x_{12},x_{13},x_{14},x_{15},x_{23},x_{24},x_{34})$. Let $\Delta_{1}=\chi_{1}(C),\; \Delta_{2}=\chi_{2}(C)\ \Delta_{1}'=\chi_{1}(C'),\ \mbox{and}\ \Delta_{2}'=\chi_{2}(C').$ We present $TLb_{8}\ $ as a union of the following disjoint subsets:
 
 $U_{8}^{1}=\{L(C)\in TLb_{8}\ :c_{34}\neq 0,\ c_{11} \neq 0,\ c_{12}\neq 0\};$
 
 $U_{8}^{2}=\{L(C)\in TLb_{8}\ : c_{34}\neq 0,\ c_{11}\neq 0,\
 c_{12}=0,\Delta_{1}\neq 0, c_{23}\neq 0\};$
 
 $U_{8}^{3}=\{L(C)\in TLb_{8}\ : c_{34}\neq 0,\
 c_{11}\neq 0,c_{12}=0,\Delta_{1}\neq 0,c_{23}=0,c_{13}\neq 0\};$
 
 $U_{8}^{4}=\{L(C)\in TLb_{8}\ : c_{34}\neq 0,\
  c_{11}\neq 0,c_{12}=0,\Delta_{1}\neq 0,c_{23}=0,c_{13}=0\};$

 $U_{8}^{5}=\{L(C)\in TLb_{8}\ : c_{34}\neq 0,\ c_{11}\neq 0,\
  c_{12}=0,\Delta_{1}=0, c_{23}\neq 0\};$
 
 $U_{8}^{6}=\{L(C)\in TLb_{8}\ : c_{34}\neq 0,\
  c_{11}\neq 0,c_{12}=0,\Delta_{1}=0,c_{23}=0,c_{13}\neq 0\};$
  
 $U_{8}^{7}=\{L(C)\in TLb_{8}\ : c_{34}\neq 0,\
   c_{11}\neq 0,c_{12}=0,\Delta_{1}=0,c_{23}=0,c_{13}=0\};$
 
 $U_{8}^{8}=\{L(C)\in TLb_{8}\ :c_{34}\neq 0,\ c_{11}=0,c_{01}\neq 0,c_{12}\neq 0\};$
 
 $U_{8}^{9}=\{L(C)\in TLb_{8}\ :c_{34}\neq 0,\ c_{11}=0,c_{01}\neq 0,c_{12}=0,c_{23}\neq 0\};$
 
 $U_{8}^{10}=\{L(C)\in TLb_{8}\ :c_{34}\neq 0,\ c_{11}=0,c_{01}\neq 0,c_{12}=0,c_{23}=0,c_{13}\neq 0\};$
 
 $U_{8}^{11}=\{L(C)\in TLb_{8}\ :c_{34}\neq 0,\ c_{11}=0,c_{01}\neq 0,c_{12}=0,c_{23}=0,c_{13}=0\};$
 
 $U_{8}^{12}=\{L(C)\in TLb_{8}\ :c_{34}\neq 0,\ c_{11}=0,c_{01}=0,c_{00}\neq 0,c_{12}\neq 0\};$
 
 $U_{8}^{13}=\{L(C)\in TLb_{8}\ :c_{34}\neq 0,\ c_{11}=c_{01}=0,c_{00}\neq 0,c_{12}=0,c_{13}\neq 0,c_{23}\neq 0\};$
 
 $U_{8}^{14}=\{L(C)\in TLb_{8}\ :c_{34}\neq 0,\ c_{11}=c_{01}=0,c_{00}\neq 0,c_{12}=0,c_{13}\neq 0,c_{23}=0\};$
 
 $U_{8}^{15}=\{L(C)\in TLb_{8}\ :c_{34}\neq 0,\ c_{11}=c_{01}=0,c_{00}\neq 0,c_{12}=c_{13}=0,c_{23}\neq 0\};$

 $U_{8}^{16}=\{L(C)\in TLb_{8}\ :c_{34}\neq 0,\ c_{11}=0,c_{01}=0,c_{00}\neq 0,c_{12}=0,c_{13}=0,c_{23}=0,c_{14}\neq 0\};$
 
 $U_{8}^{17}=\{L(C)\in TLb_{8}\ :c_{34}\neq 0,\ c_{11}=0,c_{01}=0,c_{00}\neq 0,c_{12}=0,c_{13}=0,c_{23}=0,c_{14}=0\};$
 
 $U_{8}^{18}=\{L(C)\in TLb_{8}\ :c_{34}\neq 0,\ c_{11}=0,c_{01}=0,c_{00}=0,c_{12}\neq 0,c_{13}\neq 0\};$
 
 $U_{8}^{19}=\{L(C)\in TLb_{8}\ :c_{34}\neq 0,\ c_{11}=0,c_{01}=0,c_{00}=0,c_{12}\neq 0,c_{13}=0\};$
 
 
 $U_{8}^{20}=\{L(C)\in TLb_{8}\
 :c_{34}\neq 0,\ c_{11}=c_{01}=c_{00}=c_{12}=0,c_{13}\neq 0,c_{23}\neq 0\};$
 
 $U_{8}^{21}=\{L(C)\in TLb_{8}\
 :c_{34}\neq 0,\ c_{11}=c_{01}=c_{00}=c_{12}=0,c_{13}\neq 0,c_{23}=0\};$
 
 $U_{8}^{22}=\{L(C)\in TLb_{8}\
 :c_{34}\neq 0,\ c_{11}=c_{01}=c_{00}=c_{12}=c_{13}=0,c_{23}\neq 0\};$
 
 $U_{8}^{23}=\{L(C)\in TLb_{8}\
 :c_{34}\neq 0,\ c_{11}=c_{01}=c_{00}=c_{12}=c_{13}=c_{23}=0,c_{14}\neq 0\};$
 
 $U_{8}^{24}=\{L(C)\in TLb_{8}\
 :c_{34}\neq 0,\ c_{11}=c_{01}=c_{00}=c_{12}=c_{13}=c_{23}=c_{14}=0\};$
 
 $U_{8}^{25}=\{L(C)\in TLb_{8}\
 :c_{34}=0,c_{23}\neq 0,c_{11}\neq 0,c_{12}\neq 0\}.$
 
 $U_{8}^{26}=\{L(C)\in TLb_{8}\
 :c_{34}=0,c_{23}\neq 0,c_{11}\neq 0,c_{12}=0,\Delta_{1}\neq 0\}.$
 
 $U_{8}^{27}=\{L(C)\in TLb_{8}\
 :c_{34}=0,c_{23}\neq 0,c_{11}\neq 0,c_{12}=0,\Delta_{1}=0\}.$
 
 $U_{8}^{28}=\{L(C)\in TLb_{8}\
 :c_{34}=0,c_{23}\neq 0,c_{11}=0,c_{01}\neq 0,c_{12}\neq 0\}.$
 
 $U_{8}^{29}=\{L(C)\in TLb_{8}\
 :c_{34}=0,c_{23}\neq 0,c_{11}=0,c_{01}\neq 0,c_{12}=0\}.$
 
 $U_{8}^{30}=\{L(C)\in TLb_{8}\ :c_{34}=0,c_{23}\neq 0,c_{11}=0,c_{01}=0,c_{00}\neq 0,c_{12}\neq 0\};$
 
 $U_{8}^{31}=\{L(C)\in TLb_{8}\ :c_{34}=0,c_{23}\neq 0,c_{11}=0,c_{01}=0,c_{00}\neq 0,c_{12}=0\};$
 
 $U_{8}^{32}=\{L(C)\in TLb_{8}\ :c_{34}=0,c_{23}\neq 0,c_{11}=c_{01}=c_{00}=0,c_{12}\neq 0,c_{13}\neq 0\};$
 
 $U_{8}^{33}=\{L(C)\in TLb_{8}\ : c_{34}=0,c_{23}\neq 0,c_{11}=c_{01}=c_{00}=0,c_{12}\neq 0,c_{13}=0\};$

 $U_{8}^{34}=\{L(C)\in TLb_{8}\ :c_{34}=0,c_{23}\neq 0,c_{11}=c_{01}=c_{00}=c_{12}=0,c_{13}\neq 0\};$
 
 $U_{8}^{35}=\{L(C)\in TLb_{8}\ :c_{34}=0,c_{23}\neq 0,c_{11}=c_{01}=c_{00}=c_{12}=c_{13}=0\};$
 
 $U_{8}^{36}=\{L(C)\in TLb_{8}\ :c_{34}=0,c_{23}=0,c_{24}\neq 0,c_{11}\neq 0,c_{12}\neq 0\};$
 
 $U_{8}^{37}=\{L(C)\in TLb_{8}\ :c_{34}=0,c_{23}=0,c_{24}\neq 0,c_{11}\neq 0,c_{12}=0,\Delta_{1}\neq 0,c_{13}\neq 0\};$
 
  $U_{8}^{38}=\{L(C)\in TLb_{8}\ :c_{34}=0,c_{23}=0,c_{24}\neq 0,c_{11}\neq 0,c_{12}=0,\Delta_{1}\neq 0,c_{13}=0\};$
 
 $U_{8}^{39}=\{L(C)\in TLb_{8}\ :c_{34}=0,c_{23}=0,c_{24}\neq 0,c_{11}\neq 0,c_{12}=0,\Delta_{1}=0,c_{13}\neq 0\};$
 
 $U_{8}^{40}=\{L(C)\in TLb_{8}\ : c_{34}=0,c_{23}=0,c_{24}\neq 0,c_{11}\neq 0,c_{12}=0,\Delta_{1}=0,c_{13}=0\};$
 
 $U_{8}^{41}=\{L(C)\in TLb_{8}\ :c_{34}=0,c_{23}=0,c_{24}\neq 0,c_{11}=0,c_{01}\neq 0,c_{12}\neq 0\};$
 
 $U_{8}^{42}=\{L(C)\in TLb_{8}\ :c_{34}=c_{23}=0,c_{24}\neq 0,c_{11}=0,c_{01}\neq 0,c_{12}=0,c_{13}\neq 0\};$
 
 $U_{8}^{43}=\{L(C)\in TLb_{8}\ :c_{34}=c_{23}=0,c_{24}\neq 0,c_{11}=0,c_{01}\neq 0,c_{12}=c_{13}=0\};$
 
 $U_{8}^{44}=\{L(C)\in TLb_{8}\ :c_{34}=c_{23}=0,c_{24}\neq 0,c_{11}=c_{01}=0,c_{00}\neq 0,c_{12}\neq 0\};$
 
 $U_{8}^{45}=\{L(C)\in TLb_{8}\ :c_{34}=0,c_{23}=0,c_{24}\neq 0,c_{11}=0,c_{01}=0,c_{00}\neq 0, c_{12}=0,c_{13}\neq 0\};$
 
 $U_{8}^{46}=\{L(C)\in TLb_{8}\ :c_{34}=0,c_{23}=0,c_{24}\neq 0,c_{11}=0,c_{01}=0,c_{00}\neq 0,c_{12}=0,c_{13}=0\};$
 
 $U_{8}^{47}=\{L(C)\in TLb_{8}\ :c_{34}=c_{23}=0,c_{24}\neq 0,c_{11}=c_{01}=c_{00}=0,c_{12}\neq 0\};$
 
 $U_{8}^{48}=\{L(C)\in TLb_{8}\ :c_{34}=0,c_{23}=0,c_{24}\neq 0,c_{11}=0,c_{01}=0,c_{00}=0,c_{12}=0,c_{13}\neq 0,c_{14}\neq 0\};$
 
 $U_{8}^{49}=\{L(C)\in TLb_{8}\ :c_{34}=0,c_{23}=0,c_{24}\neq 0,c_{11}=0,c_{01}=0,c_{00}=0,c_{12}=0,c_{13}\neq 0,c_{14}=0\};$
 
 $U_{8}^{50}=\{L(C)\in TLb_{8}\
 :c_{34}=c_{23}=0,c_{24}\neq 0,c_{11}=c_{01}=c_{00}=c_{12}=c_{13}=0\};$
 
 $U_{8}^{51}=\{L(C)\in TLb_{8}\
 :c_{34}=c_{23}=c_{24}=0,c_{11}\neq 0,c_{12}\neq 0\};$
 
 $U_{8}^{52}=\{L(C)\in TLb_{8}\
 :c_{34}=c_{23}=c_{24}=0,c_{11}\neq 0,c_{12}=0,\Delta_{1}\neq 0,c_{13}\neq 0\};$
 
 $U_{8}^{53}=\{L(C)\in TLb_{8}\
 :c_{34}=0,c_{23}=0,c_{24}=0,c_{11}\neq 0,c_{12}=0,\Delta_{1}\neq 0,c_{13}=0,c_{14}\neq 0\};$
 
 $U_{8}^{54}=\{L(C)\in TLb_{8}\
 :c_{34}=0,c_{23}=0,c_{24}=0,c_{11}\neq 0,c_{12}=0,\Delta_{1}\neq 0,c_{13}=0,c_{14}=0\};$
 
 $U_{8}^{55}=\{L(C)\in TLb_{8}\
 :c_{34}=c_{23}=c_{24}=0,c_{11}\neq 0,c_{12}=0,\Delta_{1}=0,c_{13}\neq 0\};$
 
 $U_{8}^{56}=\{L(C)\in TLb_{8}\
 :c_{34}=c_{23}=c_{24}=0,c_{11}\neq 0,c_{12}=\Delta_{1}=c_{13}=0,c_{14}\neq 0\};$
 
 $U_{8}^{57}=\{L(C)\in TLb_{8}\
 :c_{34}=0,c_{23}=0,c_{24}=0,c_{11}\neq 0,c_{12}=0,\Delta_{1}=0,c_{13}=0,c_{14}=0,c_{15}\neq 0\};$
 
 $U_{8}^{58}=\{L(C)\in TLb_{8}\
 :c_{34}=0,c_{23}=0,c_{24}=0,c_{11}\neq 0,c_{12}=0,\Delta_{1}=0,c_{13}=0,c_{14}=0,c_{15}=0\};$
 
 $U_{8}^{59}=\{L(C)\in TLb_{8}\
 :c_{34}=c_{23}=c_{24}=c_{11}=0,c_{01}\neq 0,c_{12}\neq 0\};$
 
 $U_{8}^{60}=\{L(C)\in TLb_{8}\
 :c_{34}=c_{23}=c_{24}=c_{11}=0,c_{01}\neq 0,c_{12}=0,c_{13}\neq 0\};$
 
 $U_{8}^{61}=\{L(C)\in TLb_{8}\ :c_{34}=0,c_{23}=0,c_{24}=0,c_{11}=0,c_{01}\neq 0,c_{12}=0,c_{13}=0,c_{14}\neq 0\};$
 
 $U_{8}^{62}=\{L(C)\in TLb_{8}\ :c_{34}=0,c_{23}=0,c_{24}=0,c_{11}=0,c_{01}\neq 0,c_{12}=0,c_{13}=0,c_{14}=0,c_{15}\neq 0\};$
 
 $U_{8}^{63}=\{L(C)\in TLb_{8}\ : c_{34}=0,c_{23}=0,c_{24}=0,c_{11}=0,c_{01}\neq0,c_{12}=0,c_{13}=0,c_{14}=0,c_{15}=0\};$
 
 $U_{7}^{64}=\{L(C)\in TLb_{7}\ :c_{34}=c_{23}=c_{24}=c_{11}=c_{01}=0,c_{00}\neq 0,c_{12}\neq 0\};$ 
 
 $U_{8}^{65}=\{L(C)\in TLb_{8}\ : c_{34}=0,c_{23}=0,c_{24}=0,c_{11}=0,c_{01}=0,c_{00}\neq 0,c_{12}=0,c_{13}\neq 0\};$
 
 $U_{8}^{66}=\{L(C)\in TLb_{8}\ :c_{34}=0,c_{23}=0,c_{24}=0,c_{11}=0,c_{01}=0,c_{00}\neq 0,c_{12}=0,c_{13}=0,c_{14}\neq 0\};$
 
 $U_{8}^{67}=\{L(C)\in TLb_{8}\ : c_{34}=0,c_{23}=0,c_{24}=0,c_{11}=0,c_{01}=0,c_{00}\neq 0,$ 
  
  \hspace{7.5cm}$\left.c_{12}=0,c_{13}=0,c_{14}=0,c_{15}\neq 0\right\};$
 
 $U_{8}^{68}=\{L(C)\in TLb_{8}\ :c_{34}=0,c_{23}=0,c_{24}=0,c_{11}=0,c_{01}=0,c_{00}\neq 0,$ 
  
  \hspace{8cm}$\left.c_{12}=0,c_{13}=c_{14}=0,c_{15}=0\right\};$
 
 $U_{8}^{69}=\{L(C)\in TLb_{8}\ :c_{34}=c_{23}=c_{24}=c_{11}=c_{01}=c_{00}=0,c_{12}\neq 0\};$
 
 $U_{8}^{70}=\{L(C)\in TLb_{8}\ :c_{34}=c_{23}=c_{24}=c_{11}=c_{01}=c_{00}=c_{12}=0,c_{13}\neq 0\};$

 $U_{8}^{71}=\{L(C)\in TLb_{8}\ :c_{34}=0,c_{23}=0,c_{24}=0,c_{11}=0,c_{01}=0,c_{00}=0,$ 
  
  \hspace{9.5cm}$\left.c_{12}=0,c_{13}=0,c_{14}\neq 0\right\};$
 
 $U_{8}^{72}=\{L(C)\in TLb_{8}\ :c_{34}=0,c_{23}=0,c_{24}=0,c_{11}=0,c_{01}=0,c_{00}=0 $ 
 
 \hspace{8cm}$\left. c_{12}=0,c_{13}=0,c_{14}=0,c_{15}\neq 0\right\};$
 
 $U_{8}^{73}=\left\{L(C)\in TLb_{8}\ :c_{34}=c_{23}=c_{24}=c_{11}=c_{01}=c_{00}=c_{12}=c_{13}=c_{14}=c_{15}=0\right\}.$\\

 The above subsets are disjoint. Some of them are single orbits such as $U_{8}^{14},$ $U_{8}^{15},$ $U_{8}^{16},$ $U_{8}^{17},$  $U_{8}^{21},$ $U_{8}^{22},$ $U_{8}^{23},$ $U_{8}^{24},$ $U_{8}^{25},$ $U_{8}^{35}$ $U_{8}^{36},$ $U_{8}^{48},$ $U_{8}^{52},$ $U_{8}^{57}$ while some are parametric family of orbits like $U_{8}^{1},$ $U_{8}^{2},$ $U_{8}^{3},$  $U_{8}^{4},$ $U_{8}^{5},$ $U_{8}^{7},$ $U_{8}^{8},$  $U_{8}^{10},$ $U_{8}^{13},$ $U_{8}^{18},$ among others.\\ The following proposition gives the isomorphism criterion only for one subset from each type, that is one from infinitely many orbits case and one from single orbits case.
 \newpage
 \begin{pro}\emph{}
 \begin{enumerate}\item[\em{1.}]
 \begin{enumerate}
 \item[\em{(a)}]  Two algebras $L(C')$ and $L(C)$ from $U_{8}^{6}$ are
 isomorphic, if and only if \\ $$ \frac{c'^{5}_{13}\Delta_{2}'^{5}}{c'^{7}_{11}c'^{3}_{34}}=\frac{c^{5}_{13}\Delta_{2}^{5}}{c_{11}^{7}c_{34}^{3}}$$
  \item[\em{(b)}] For any $\lambda \in \mathbb{C}, $  there exists $L(C)\in
 U_{8}^{6}:$ 
 $$\frac{c^{5}_{13}\Delta_{2}^{5}}{c_{11}^{7}c_{34}^{3}}=\lambda.$$ {\em Then orbits in
 $U_{8}^{6}$ can be parameterized as} $L\left(0,0,1,0,\lambda,0,0,0,0,1\right), \ \lambda \in \mathbb{C}.$
 \end{enumerate}
 \item[\em{2.}] The subset $ U_{8}^{14}$ with the representative \ $ L(1,0,0,0,1,0,0,0,0,1)$ is a single orbit.
 \end{enumerate}
 \begin{proof}
 \begin{enumerate}\item
 \begin{enumerate}
 \item\ ``If'' part is due to {\bf Theorem \ref{T3}} if one substitutes the
 expressions for $c_{01}'$ $c_{11}',$ $c_{13}',$ \ and $c_{34}'$ into
 $$ \frac{c'^{5}_{13}\Delta_{2}'^{5}}{c'^{7}_{11}c'^{3}_{34}}.$$
 
`` Only if '' part. \ Let the equalities $$ \frac{c'^{5}_{13}\Delta_{2}'^{5}}{c'^{7}_{11}c'^{3}_{34}}=\frac{c^{5}_{13}\Delta_{2}^{5}}{c_{11}^{7}c_{34}^{3}}$$
 holds.
 
Consider the base change (\ref{BC}) in the proof of {\bf Theorem \ref{T3}} with
 $$ A_{0}^{5}=\frac{c_{11}}{c_{34}}, A_1={\frac {-A_{{0}}c_{01}}{2\,c_{11}}}, \
 B^{5}_1=\frac{\left(2c_{11}-c_{01}c_{34}\right)^{5}}{2c_{11}^{4}c_{34}^{6}},$$ and 
 \begin{align*}
  B_3=&\frac{A^{3}_{0}B^{2}_{2}c_{34}+A^{2}_{0}A_{1}B^{2}_{1}c_{13}c_{34}-A^{3}_{0}B_{1}^{2}c_{24}}{2A^{3}_{0}B_1c_{34}}
 \end{align*}

 This changing transforms $L(C)$ into
 $$ L\left(0, 0, 1, 0, \frac{c^{5}_{13}\Delta_{2}^{5}}{c_{11}^{7}c_{34}^{3}},0,0,0,0,1\right).$$ An analogous base change transforms $L(C')$ into $$ L\left(0, 0, 1, 0, \frac{c'^{5}_{13}\Delta_{2}'^{5}}{c'^{7}_{11}c'^{3}_{34}},0,0,0,0,1\right).$$ Since
 $$\frac{c'^{5}_{13}\Delta_{2}'^{5}}{c'^{7}_{11}c'^{3}_{34}}=\frac{c^{5}_{13}\Delta_{2}^{5}}{c_{11}^{7}c_{34}^{3}}$$ the algebras are isomorphic. 
 \item\ It is an obvious consequence of 1(a). 
  \end{enumerate}
  \item This is proved by putting appropriate values of $A_0$, $A_1$, $B_1$, $B_2$, and $B_3$ into the base change (\ref{BC}). It is not too laborious to realise that $A^{9}_{0}=\frac{c_{00}c^{2}_{13}}{c_{34}}$,\ $B_{1}=\frac{A^{3}_{0}}{c_{13}}$,\ $B_{2}=\frac{B_{1}c_{14}}{c_{13}}$,\ $B_{3}=\frac{A^{3}_{0}B^{2}_{2}c_{34}-A^{3}_{0}B^{2}_{1}c_{24}}{2A^{3}_{0}B_{1}c_{34}}$,\\and $A_{1},\ B_{4},\ B_{5}\in \mathbb{C}$.  
  \end{enumerate}
 \end{proof}
 \end{pro}
\section{Conclusion}
This section summarizes the result of this study by listing its milestones. The isomorphism classes in $TLb_7$ and $TLb_8$ are presented, together with the invariant functions in parametric family of orbits case, in four tables as shown in Table \ref{tab:1}, Table \ref{tab:2}, Table \ref{tab:3} and Table \ref{tab:4} below.
\begin{enumerate}
\item In $TLb_7$ we outline $30$ isomorphism classes ($10$ parametric family of orbits and $20$ single orbits) of seven dimensional complex filiform Leibniz algebras while in $TLb_8$, there $73$ isomorphism classes ($32$ single orbits and $41$ parametric family of orbits). These classes exhaust all possible cases.
\item The seven-dimensional filiform Lie algebras are covered by $U_{7}^{10},$
$U_{7}^{11},$ $U_{7}^{12},$
$U_{7}^{26},$ $U_{7}^{27},$ $U_{7}^{28},$ $U_{7}^{29},$ and $U_{7}^{30}$, while the eight-dimensional filiform Lie algebras are covered by $U_{8}^{18},$ $U_{8}^{19},$ $U_{8}^{20},$ $U_{8}^{21},$ $U_{8}^{22},$ $U_{8}^{23},$ $U_{8}^{24},$ $U_{8}^{32},$ $U_{8}^{33},$ $U_{8}^{34},$ $U_{8}^{35},$ $U_{8}^{47},$ $U_{8}^{48},$ $U_{8}^{49},$ $U_{8}^{50},$ $U_{8}^{69},$ $U_{8}^{70},$ $U_{8}^{71},$ $U_{8}^{72}$ and $U_{8}^{73}.$
\item The list of resulting filiform Lie algebras in seven-dimensional case agrees with that given in  \cite{AG} and \cite{GJKh}. However, we found that the list in $TLb_8$ agrees with the list given in \cite{AG} while we have one more class of filiform Lie algebras than the classes given in \cite{GJKh}. Therefore, the list of isomorphism classes given in \cite{GJKh} should be accordingly corrected.
\end{enumerate}

\begin{table}[h!]
\caption{The list of single orbits in $TLb_7$}
\label{tab:1} 
\begin{tabular}{p{2cm}p{4.9cm}p{2cm}p{4.9cm}}\hline
  
  Subsets &  Representative of orbits & Subsets &  Representative of orbits  \\
  \hline\\
   $U_{7}^{6} \ $ & $ L(0,1,0,0,0,0,1)$ & $U_{7}^{9} \ $ & $ L(1,0,0,0,0,0,1)$ \\  
   $U_{7}^{11} $ & $ L(0,0,0,0,1,0,1)$ & $U_{7}^{21} $ & $ L(0,1,0,0,0,0,0)$ \\  
   $U_{7}^{12} $& $L(0,0,0,0,0,0,1) $ & $U_{7}^{22} $ & $  L(1,0,0,1,0,0,0)$ \\ 
    $U_{7}^{15} $ & $ L(0,0,1,0,1,0,0)$ & $U_{7}^{23} $ & $ L(1,0,0,0,1,0,0)$ \\ 
   $U_{7}^{16} $ & $L(0,0,1,0,0,1,0) $ & $U_{7}^{24} $ & $ L(1,0,0,0,0,1,0)$ \\
   $U_{7}^{17} $ & $  L(0,0,1,0,0,0,0) $ & $U_{7}^{25} $ & $L(1,0,0,0,0,0,0)$ \\
   $U_{7}^{18} $ & $  L(0,1,0,1,0,0,0) $ & $U_{7}^{26}$ & $ L(0,0,0,1,0,0,0) $ \\ 
   $U_{7}^{19} $ & $  L(0,1,0,0,1,0,0) $ &  $U_{7}^{27}$ & $ L(0,0,0,0,1,0,0) $ \\ 
   $U_{7}^{20} $ & $  L(0,1,0,0,0,1,0) $ &  $U_{7}^{28}$ & $ L(0,0,0,0,1,1,0) $ \\ 
   $U_{7}^{30} $ & $  L(0,0,0,0,0,0,0) $ &  $U_{7}^{29}$ & $ L(0,0,0,0,0,1,0) $ \\ \hline
\end{tabular}
\end{table}
\begin{table}[h!]
\caption{Parametric family of orbits of $TLb_7$ with invariant functions}
\label{tab:2}
\begin{tabular}{p{2cm}p{4cm}p{8cm}}
  \hline
  Subsets &  Representative of orbits& Invariant functions\\
  \hline
  $U_{7}^{1}$ & $L\left(\lambda_1, 0, 1, \lambda_2,0, 0,1\right)$ & $ f_1(X)= \left(\frac
{x_{23}}{x_{11}}\right)^{8}\,\left(\chi_{1}(X)\right)^{3},$
$  f_2(X)={\frac {x_{12}}{x_{23}}}$  \\ 

  $U_{7}^{2}$ & $ L\left(1, 0, 1,0,\lambda_1,0,\lambda_2\right)$ & $ f_1(X)=\left(\frac {x_{13}}{x_{11}}\right)^{4} \left(\chi_{1}(X)\right),$
$f_2(X)=\left(\frac{x_{23}}{x_{11}}\right)^{8}\left(\chi_{1}(X)\right)^{3}$ \\ 

  $U_{7}^{3}$ & $L\left(0, 0,1,0, \lambda, 0,1\right)$ & $f(X)={\frac{x_{13}^{3}}{{x_{11}}{x_{23}^{2}}}}$  \\

  $U_{7}^{4}$ & $L\left(0,1,0, \lambda, 0, 0, 1\right) $ & $ f(X)={\frac {x_{12}}{x_{23}}}$  \\
  
  $U_{7}^{5}$ & $  L\left(0,\lambda, 0, 0,1, 0, 1\right)$ & $ f(X)=x_{01}\left({\frac
{x_{23}}{{x_{13}}}}\right)^{{4}} $ \\
 
   $U_{7}^{7}$ & $
   L\left(1,0,0,1,0,0,\lambda\right)$ & $
   f(X)=\frac{x_{23}}{x_{12}}  $ \\

   $U_{7}^{8}$ & $
L\left(\lambda,0,0,0,1,0,1\right)$ & $
f(X)=\frac{x_{00}{x_{23}^{6}}}{x_{13}^{5}}  $  \\ 
 $U_{7}^{10}$ & $L\left(0,0,0 ,1,0,0,\lambda\right) $ & $ f(X)=\frac {x_{23}}{x_{12}}$\\ 
     $U_{7}^{13}$ & $ L\left(\lambda,0,1,1,0,0,0\right)$ & $
f(X)=\left(\frac
{x_{12}}{x_{11}}\right)^{8}\,\left(\chi_{1}(X)\right)^{3}$\\

   $U_{7}^{14}$ & $ L\left(1,0,1,0,\lambda,0,0,\right) $ &
$f(X)=\left(\frac{x_{13}}{x_{11}}\right)^{4}\left(\chi_{1}(X)\right)$\\ \hline
\end{tabular}
\end{table}
\qquad \qquad \qquad \qquad \qquad \qquad \qquad \qquad \qquad $X=(x_{00},x_{01},x_{11},x_{12},x_{13},x_{14},x_{23})$.
\begin{table}[h!]
\caption{The list of single orbits in $TLb_{8}$}
\label{tab:3} 
\begin{tabular}{p{2cm}p{4.9cm}p{2cm}p{4.9cm}}\hline
  
  Subsets &  Representative of orbits & Subsets &  Representative of orbits  \\
  \hline\\
   $U_{8}^{14} \ $ & $ L(1,0,0,0,1,0,0,0,0,1)$ & $U_{8}^{15} \ $ & $ L(1,0,0,0,0,0,0,1,0,1)$ \\[2ex]  
   $U_{8}^{16} $ & $ L(1,0,0,0,0,1,0,0,0,1)$ & $U_{8}^{17} $ & $ L(1,0,0,0,0,0,0,0,0,1)$ \\ [2ex] 
   $U_{8}^{20} $& $L(0,0,0,0,1,0,0,1,0,1)$ & $U_{8}^{21} $ & $ L(0,0,0,0,1,0,0,0,0,1)$ \\[2ex] 
    $U_{8}^{22} $ & $ L(0,0,0,0,0,0,0,1,0,1)$ & $U_{8}^{23} $ & $ L(0,0,0,0,0,1,0,0,0,1)$ \\[2ex] 
   $U_{8}^{24} $ & $L(0,0,0,0,0,0,0,0,0,1) $ & $U_{8}^{34} $ & $ L(0,0,0,0,1,0,0,1,0,0)$ \\[2ex]
   $U_{8}^{35} $ & $ L(0,0,0,0,0,0,0,1,0,0)$ & $U_{8}^{47} $ & $ L(0,0,0,1,0,0,0,0,1,0)$ \\[2ex]
   $U_{8}^{50} $ & $ L(0,0,0,0,0,1,0,0,1,0)$ & $U_{8}^{55}$ & $ L(0,0,1,0,1,0,0,0,0,0)$ \\ [2ex]
   $U_{8}^{56} $ & $ L(0,0,1,0,0,1,0,0,0,0) $ &  $U_{8}^{57}$ & $ L(0,0,1,0,0,0,1,0,0,0) $ \\[2ex] 
   $U_{8}^{58} $ & $ L(0,0,1,0,0,0,0,0,0,0) $ &  $U_{8}^{59}$ & $ L(0,1,0,1,0,0,0,0,0,0) $ \\[2ex] 
   $U_{8}^{60} $ & $ L(0,1,0,0,1,0,0,0,0,0) $ &  $U_{8}^{61}$ & $ L(0,1,0,0,0,1,0,0,0,0) $ \\[2ex]
   $U_{8}^{62} $ & $ L(0,1,0,0,0,0,1,0,0,0) $ &  $U_{8}^{63}$ & $ L(0,1,0,0,0,0,0,0,0,0) $ \\[2ex] $U_{8}^{64} $ & $ L(1,0,0,1,0,0,0,0,0,0) $ &  $U_{8}^{65}$ & $ L(1,0,0,0,1,0,0,0,0,0) $ \\[2ex] $U_{8}^{66} $ & $ L(1,0,0,0,0,1,0,0,0,0) $ &  $U_{8}^{67}$ & $ L(1,0,0,0,0,0,1,0,0,0) $ \\[2ex]
   $U_{8}^{68} $ & $ L(1,0,0,0,0,0,0,0,0,0) $ &  $U_{8}^{69}$ & $ L(0,0,0,1,0,0,0,0,0,0) $ \\[2ex]
   $U_{8}^{70} $ & $ L(0,0,0,0,1,0,0,0,0,0) $ &  $U_{8}^{71}$ & $ L(0,0,0,0,0,1,0,0,0,0) $ \\[2ex]
   $U_{8}^{72} $ & $ L(0,0,0,0,0,0,1,0,0,0) $ &  $U_{8}^{73}$ & $ L(0,0,0,0,0,0,0,0,0,0) $ \\ \hline
\end{tabular}
\end{table}

\begin{table}[h!]
\caption{Parametric family of orbits in $TLb_8$ with invariant functions}
\label{tab:4}
\begin{tabular}{p{3cm}p{6cm}p{6.5cm}}
  \hline
  Subsets &  Representative of orbits & Invariant functions \\
  \hline
 $U_{8}^{1}$ & ${L\left(\lambda_1, 0, 1, \lambda_2,0,0,0,\lambda_3,\lambda_4,1\right)}$ & $ \begin{aligned}f_1(X)&= \left(\frac{x_{34}}{\chi_{2}(X)}\right)^{2}\,\chi_{1}(X),\\
f_2(X)&=\frac{x^{5}_{12}(\chi_{2}(X))^{5}}{x_{11}^{6}x_{34}^{4}},\\
f_3(X)&=\frac{x^{5}_{23}(\chi_{2}(X))^{5}}{x_{11}^{6}x_{34}^{4}},\\
f_4(X)&=\frac{(2x_{11}x_{24}+x_{01}(\chi_{4}(X)))^{5}}{2x^{11}_{11}x^{9}_{34}}\end{aligned}$\\ \hline

 $U_{8}^{2}$ & $ L\left(1, 0, 1,0,\lambda_1,0,0,\lambda_2,0,\lambda_3\right)$ & $\begin{aligned} f_1(X)&=\frac {x^{10}_{13}(\chi_{2}(X))^{10}}{x^{14}_{11}}\left(\chi_{1}(x)\right)^{3},\\
f_2(X)&=\left(\frac{x_{23}}{x_{11}}\right)^{2}\left(\chi_{1}(x)\right),\\ f_3(X)&= \left(\frac{x_{34}}{\chi_{2}(X)}\right)^{2} \,\left(\chi_{1}(x)\right),\end{aligned}$
\\ \hline
  \end{tabular}
    \end{table}
  \clearpage
    \begin{table}[h!]
    \ContinuedFloat
    \begin{tabular}{p{3cm}p{6cm}p{6.5cm}}
      \hline
      Subsets &  Representative of orbits& Invariant functions\\
      \hline

  $U_{8}^{3}$ & $L\left(1, 0, 1,0,\lambda_1,0,0,0,0,\lambda_2\right)$ & $\begin{aligned}
  f_1(X)&=\frac{x^{10}_{13}(\chi_{2}(X))^{4}}{x^{14}_{11}} \left(\chi_{1}(X)\right)^{3},\\
  f_{2}(X)&=\left(\frac{x_{34}}{\chi_{2}(X)}\right)^{2}\,\left(\chi_{1}(X)\right)\end{aligned}$\\ \hline

  $U_{8}^{4}$ & $L\left(1, 0, 1,0,0,\lambda_1,0,0,0,\lambda_2\right)$ & 
  $\begin{aligned}f_1(X)&=\frac{x^{5}_{14}(\chi_{2}(X))^{3}}{x^{8}_{11}} \left(\chi_{1}(X)\right),\\
   f_{2}(X)&=\left(\frac{x_{34}}{\chi_{2}(X)}\right)^{2}\,\left(\chi_{1}(X)\right)\end{aligned}$  \\ \hline

    $U_{8}^{5}$ & $ L\left(0, 0, 1,0,\lambda_1,0,0,\lambda_2,0,1\right)$ & $\begin{aligned} f_1(X)&=\frac{x^{5}_{13}(\chi_{2}(X))^{5}}{x_{11}^{7}x_{34}^{3}},\\
  f_2(X)&=\frac{x^{5}_{23}(\chi_{2}(X))^{5}}{x_{11}^{6}x_{34}^{4}}\end{aligned}$\\ \hline
 
   $U_{8}^{6}$ & $
  L\left(0,0,1,0,\lambda,0,0,0,0,1\right)$ & $f(X)=\frac{x^{5}_{13}(\chi_{2}(X))^{5}}{x_{11}^{7}x_{34}^{3}}$ \\ \hline
  
     $U_{8}^{7}$ & $
L\left(0,0,1,0,0,\lambda,0,0,0,1\right)$ & $\small
f(X)=\frac{x^{5}_{14}(\chi_{2}(X))^{3}}{x^{8}_{11}}\left(\chi_{1}(X)\right)$ \\ \hline
 
 $U_{8}^{8}$ & $L\left(0,1,1,\lambda_1,0,0,0,\lambda_2,\lambda_3,1\right)$ & $\begin{aligned} f_{1}(X)&=\left(\frac{x_{12}}{x_{34}}\right)^{5}\left(\frac{(\chi_{3}(X))^{6}}{x_{01}^{7}}\right),\\ f_{2}(X)&=\left(\frac{x_{23}}{x_{34}}\right)^{5}\left(\frac{(\chi_{3}(X))^{6}}{x_{01}^{7}}\right),\\ f_{3}(X)&=\frac{(x_{01}x_{24}+x_{00}(\chi_{4}(X)))^{5}}{x^{11}_{01}x^{5}_{34}}\end{aligned}$\\ \hline

     $U_{8}^{9}$ & $L\left(0,1,0,0,\lambda_1,0,0,\lambda_2,0,1\right)$ & $\begin{aligned}
f_{1}(X)&=\left(\frac{x_{13}}{x_{34}}\right)^{5}\left(\frac{(\chi_{3}(X))^{7}}{x_{01}^{9}}\right),\\ f_{2}(X)&=\left(\frac{x_{23}}{x_{34}}\right)^{5}\left(\frac{(\chi_{3}(X))^{6}}{x_{01}^{7}}\right)\end{aligned}$\\ \hline

   $U_{8}^{10}$ & $ L\left(0,1,0,0,\lambda_1,0,0,0,0,1\right)$ &
$f(X)=\left(\frac{x_{13}}{x_{34}}\right)^{5}\big(\frac{(\chi_{3}(X))^{7}}{x_{01}^{9}}\big)$\\ \hline

$U_{8}^{11}$ & $ L\left(0,1,0,0,0,\lambda,0,0,0,1\right)$ &
$f(X)=\left(\frac{x_{14}}{x_{34}}\right)^{5}\Big(\frac{(\chi_{3}(X))^{8}}{x_{01}^{11}}\Big)$\\ \hline

$U_{8}^{12}$ & $ L\left(1,0,0,1,0,0,0,\lambda_1,\lambda_2,1\right)$ &
$\begin{aligned}f_{1}(X)&=\frac{x_{23}}{x_{12}},
f_{2}(X)=\frac{x_{24}x_{34}}{x^{2}_{12}}\end{aligned}$\\ \hline

$U_{8}^{13}$ & $ L\left(0,1,0,0,1,0,0,0,0,\lambda\right)$ &
$f(X)=\frac{c^{7}_{13}c_{34}}{c_{00}c^{9}_{23}}$\\ \hline

$U_{8}^{18}$ & $ L\left(0,0,0,1,0,0,0,\lambda,0,1\right)$ &
$\begin{aligned}f(X)&=\frac{x_{23}}{x_{12}}\end{aligned}$\\ \hline

$U_{8}^{19}$ & $ L\left(0,0,0,1,1,0,0,\lambda,0,1\right)$ &
$\begin{aligned}f(X)&=\frac{x_{23}}{x_{12}}\end{aligned}$\\ \hline


$U_{8}^{25}$ & $ L\left(\lambda_1,0,1,\lambda_2,0,0,0,0,\lambda_3,0\right)$ &
$\begin{aligned}f_{1}(X)&=\left(\frac{x_{23}}{x_{11}}\right)^{5}\left(\chi_{1}(X)\right)^{2},
f_{2}(X)=\frac{x_{12}}{x_{23}},\\
f_{3}(X)&=\frac{x^{4}_{24}}{x_{11}x^{4}_{23}}\end{aligned}$ \\ \hline

$U_{8}^{26}$ & $ L\left(1,0,1,0,\lambda_1,0,0,\lambda_2,0,0\right)$ &
$\begin{aligned}f_{1}(X)&=\left(\frac{x_{13}}{x_{11}}\right)^{10}(\chi_{1}(X))^{3},\\
f_{2}(X)&=\left(\frac{x_{23}}{x_{11}}\right)^{5}(\chi_{1}(X))^{2}\end{aligned}$\\ \hline

$U_{8}^{27}$ & $ L\left(0,0,1,0,\lambda,0,0,1,0,0\right)$ &
$\begin{aligned}f(X)=\frac{x_{13}^{4}}{x_{11}x^{3}_{23}}\end{aligned}$ \\ \hline
\end{tabular}
  \end{table}
\clearpage
  \begin{table} [h!]
  \ContinuedFloat
  \begin{tabular}{p{3cm}p{6cm}p{6.5cm}}
    \hline
    Subsets &  Representative of orbits& Invariant functions\\
    \hline
$U_{8}^{28}$ & $ L\left(0,1,0,\lambda_1,0,0,0,1,\lambda_2,0\right)$ &
$\begin{aligned}f_{1}(X)&=\frac{x_{12}}{x_{23}},\\
f_{2}(X)&=\frac{x^{5}_{24}}{x_{01}x^{5}_{23}}\end{aligned}$\\ \hline

$U_{8}^{29}$ & $ L\left(0,1,0,0,\lambda,0,0,1,0,0\right)$ &
$\begin{aligned}f(X)&=\frac{x^{5}_{13}}{x_{01}x^{5}_{23}},\end{aligned}$ \\ \hline

$U_{8}^{30}$ & $ L\left(1,0,0,\lambda,0,0,0,1,0,0\right)$ &
$\begin{aligned}f(X)&=\frac{x_{12}}{x_{23}}\end{aligned}$\\ \hline

$U_{8}^{31}$ & $ L\left(1,0,0,0,\lambda,0,0,1,0,0\right)$ &
$\begin{aligned}f(X)&=\frac{x^{6}_{13}}{x_{00}x^{7}_{23}},\end{aligned}$ \\ \hline

$U_{8}^{32}$ & $ L\left(0,0,0,\lambda,0,0,0,1,0,0\right)$ &
$f(X)=\frac{x_{12}}{x_{23}}$
\\ \hline

$U_{8}^{33}$ & $ L\left(0,0,0,1,1,0,0,\lambda,0,0\right)$ &
$\begin{aligned}f(X)=\frac{x_{23}}{x_{12}},\end{aligned}$\\[1ex] \hline

$U_{8}^{36}$ & $ L\left(\lambda_1,0,1,\lambda_2,0,0,0,0,1,0\right)$ &
$\begin{aligned}f_{1}(X)&=\left(\frac{x_{24}}{x_{11}}\right)^{10}(\chi_{1}(X))^{3},\\
f_{2}(X)&=x^{3}_{12}\left(\frac{x_{11}}{x^{4}_{24}}\right)\end{aligned}$\\ \hline

$U_{8}^{37}$ & $ L\left(1,0,1,0,\lambda_1,0,0,0,\lambda_2,0\right)$ &
$\begin{aligned}f_{1}(X)&=\left(\frac{x_{13}}{x_{11}}\right)^{3}(\chi_{1}(X))^{10},\\
f_{2}(X)&=\left(\frac{x_{24}}{x_{11}}\right)^{3}\left(\chi_{1}(X)\right)^{10}\end{aligned}$\\ \hline

$U_{8}^{38}$ & $ L\left(1,0,1,0,0,\lambda_1,0,0,\lambda_2,0\right)$ &
$\begin{aligned}f_{1}(X)&=\left(\frac{x_{14}}{x_{11}}\right)(\chi_{1}(X))^{5},\\
f_{2}(X)&=\left(\frac{x_{24}}{x_{11}}\right)^{3}(\chi_{1}(X))^{10}\end{aligned}$\\ \hline

$U_{8}^{39}$ & $ L\left(0,0,1,0,\lambda,0,0,0,1,0\right)$ &
$\begin{aligned}f(X)=\frac{x_{13}}{x_{24}}\end{aligned}$ \\ \hline

$U_{8}^{40}$ & $ L\left(0,0,1,0,0,\lambda,0,0,1,0\right)$ &
$\begin{aligned}f(X)=\frac{x^{3}_{14}}{x_{11}x^{2}_{24}}\end{aligned}$\\ \hline

$U_{8}^{41}$ & $ L\left(0,1,0,\lambda,0,0,0,0,1,0\right)$ &
$\begin{aligned}f(X)=\frac{x_{01}x^{5}_{12}}{x^{5}_{24}}\end{aligned}$\\ \hline

$U_{8}^{42}$ & $ L\left(0,1,0,0,\lambda,0,0,0,1,0\right)$ &
$\begin{aligned}f(X)=\frac{x_{13}}{x_{24}}\end{aligned}$\\ \hline

$U_{8}^{43}$ & $ L\left(0,1,0,0,0,\lambda,0,0,1,0\right)$ &
$\begin{aligned}f(X)=\frac{x^{5}_{14}}{x_{01}x^{5}_{24}}\end{aligned}$\\ \hline

$U_{8}^{44}$ & $ L\left(1,0,0,\lambda,0,0,0,0,1,0\right)$ &
$\begin{aligned}f(X)=x^{7}_{12}\left(\frac{x^{5}_{00}}{x^{2}_{24}}\right)\end{aligned}$\\ \hline

$U_{8}^{45}$ & $ L\left(1,0,0,0,\lambda,0,0,0,1,0\right)$ &
$\begin{aligned}f(X)=\frac{x_{13}}{x_{24}}\end{aligned}$\\ \hline

$U_{8}^{46}$ & $ L\left(1,0,0,0,0,\lambda,0,0,1,0\right)$ &
$\begin{aligned}f(X)=\frac{x^{7}_{14}}{x_{00}x^{8}_{24}}\end{aligned}$\\ \hline

$U_{8}^{48}$ & $ L\left(0,0,0,0,\lambda,0,0,0,1,0\right)$ &
$\begin{aligned}f(X)=\frac{x_{13}}{x_{24}}\end{aligned}$\\ \hline

$U_{8}^{49}$ & $ L\left(0,0,0,0,\lambda,1,0,0,1,0\right)$ &
$\begin{aligned}f(X)=\frac{x_{13}}{x_{24}}\end{aligned}$\\ \hline

$U_{8}^{51}$ & $ L\left(\lambda,0,1,1,0,0,0,0,0,0\right)$ &
$\begin{aligned}f(X)=\left(\frac{x_{12}}{x_{11}}\right)^{5}(\chi_{1}(X))^{2}\end{aligned}$\\ \hline

$U_{8}^{52}$ & $ L\left(1,0,1,0,\lambda,0,0,0,0,0\right)$ &
$\begin{aligned}f(X)=\left(\frac{x_{13}}{x_{11}}\right)^{10}\left(\chi_{1}(X)\right)^{3}\end{aligned}$\\ \hline

$U_{8}^{53}$ & $ L\left(1,0,1,0,0,\lambda,0,0,0,0\right)$ &
$f(X)=\left(\frac{x_{14}}{x_{11}}\right)^{5}\chi_{1}(X),$\\ \hline

$U_{8}^{54}$ & $ L\left(1,0,1,0,0,0,\lambda,0,0,0\right)$ &
$f(X)=\left(\frac{x_{15}}{x_{11}}\right)^{10}\chi_{1}(X),$\\ \hline \hspace{1cm}
\end{tabular}
$X=(x_{00},x_{01},x_{11},x_{12},x_{13},x_{14},x_{15},x_{23},x_{24},x_{34})$.
\end{table}
\section{Acknowledgements}
The first author will like to thank the Institute for Mathematical Research, Universiti Putra Malaysia, for providing the laboratory where the research was conducted. This study was supported by the research grant 05-02-12-2188RU provided by Universiti Putra Malaysia (UPM), Malaysia.

\bibliographystyle{plain}

\end{document}